\documentclass[11pt,reqno]{amsart}
\usepackage{graphicx,tikz}
\usepackage{amsfonts,amsmath,amssymb}
\usetikzlibrary{matrix,arrows}
\newcommand{\C}{\mathbb C}
\newcommand{\R}{\mathbb R}
\newcommand{\N}{\mathbb N}

\newtheorem{theorem}{Theorem}[section]
\newtheorem{lemma}[theorem]{Lemma}

\newtheorem{cor}[theorem]{Corollary}

\newtheorem{rem}[theorem]{Remark}
\newtheorem{ex}[theorem]{Example}

\begin{document}

\title {Almost holomorphic curves in real analytic hypersurfaces}

\author{ Pierre Bonneau*  and Emmanuel Mazzilli**}

\begin{abstract}  
Using the theory of exterior differential systems, we study the existence of germs of pseudo-holomorphic disk in a real analytic hypersurface locally defined in a complex manifold equipped with $J$ a real analytic almost complex structure. The integrable case in $\C^n$ with $J$ the multiplication by $i$ has been intensively studied by several authors \cite{DF}, \cite{DA1} and \cite{DA2} for example. The non integrable case is  drastically different essentially due to the following fact : in generic case, there is no $J$-invariant  objects of dimension bigger than one. This simple observation leads to the non existence of some equivalents  of Segree varieties or  ideals of holomorphic functions which play a fundamental role in the complex case. Nevertheless in the almost complex case, we adopt the exterior differential system point of view of E.Cartan developed and clarified in \cite{BCGGG}.

\end{abstract}

\maketitle

\section{introduction}
\bigskip

In this paper we study the existence of germ of pseudo-holomorphic disk in a real analytic hypersurface locally defined in a complex manifold equipped with $J$ a real analytic almost complex structure. In the integrable case, the first motivation to study  such existence in the boundary of a bounded real analytic domain of $\C^n$  was the existence of subelliptic estimates for the $ \bar\partial$-Neumann problem in the domain (see \cite{DF}). More precisely, the non-existence of germ of disk in the boundary is a sufficient condition for the subelliptic estimates for  $\bar\partial$-Neumann problem in the domain. Nevertheless, in the non integrable context, we think that such problem is a nice geometric problem linked with the existence of peudo-holomorphic foliation on real analytic hypersurface embedded in an  almost complex manifold equipped with $J$ a real analytic almost complex structure. The main difference with the integrable case is the non existence of $J$-invariant manifold of dimension bigger than one for generic structure $J$ and so, it is not possible to localize the germ of disk in the intersection of complex hypersurface defined locally and not in sense of germs. This approach is the key trick in the complex case.

Namely in this paper, we adopt  the exterior differential system for pfaffian system point of view which is, to our knowledge, original and of interest for this problem. Roughly speaking, for the pfaffian system - which is the case when writing the PDE system associated to the existence of pseudo-holomorphic disk in an hypersurface -, we have two fundamental objects for the existence of solutions of the PDE system: the torsion and the dimension of the "tableaux". The nullity of the torsion gives a necessary conditions to obtain a solution of the system and to obtain sufficient conditions, the dimension of the tableau and the dimension of its prolongation play a fundamental role. In our context, we investigate carefully the two intrinsic notions and we prove when the system is in involution in the sense of Cartan (\cite{BCGGG}).

The paper is organized as follows: in the section $2$, we recall all the materials for almost holomorphic curves and for elementary almost complex geometry used in the paper. 

In the section $3$, we compute carefully the torsion which is "like an intrinsic way to see the Levi form" and the dimensions of the tableaux. We give necessary conditions for which the PDE system associated to our problem is in involution at the first step (see corollary $3.6$ and $3.7$). In fact, we study more general system than the previous one and we compute  the  prolongations of the tableau for which this  tableau is in involution (see theorem $3.2$ for precise statements). At the end of the section $3$, we study some examples in the complex case to illustrate the abstract results of the beginning of the section (see example 3.9) and we investigate the successive torsion for the prolongation of the  more general system until that its tableau is in involution. 

In the section $4$, we give directly a necessary and sufficient condition to obtain an ordinary integral element of the exterior differential system and so, by Cartan-$K\ddot{a}hler$ theorem, the existence of germ of disk in the hypersurface (see \cite{BCGGG} pp. 81-86 and theorem $4.1$, section $4$).

In the section $5$, we explain how to construct the PDE system on finite union of manifold  "quasi-equivalent" to the system of existence of germ of disk, for  which the tableau is in involution on each manifold and free torsion (see theorem $5.5$ for precise statements).

\section{Almost holomorphic curves}
We first provide $\mathbb{R}^{2}$ and $\mathbb{R}^{2n}$ with almost complex structures.\\
On an open set $D$ of $\mathbb{R}^{2},$ we define an almost complex structure $J_{0},$ that is to say, for all $x=(x_{1}, x_{2})\in D,$ we have an isomorphism $J_{0}(x)=J_{0}$ of (the tangent space to) $\mathbb{R}^{2},$ which depends $\mathbb{R}$-analytically to $x$ and verifies $J_{0}^{2}=-Id$ ($Id$ is the identity). We note $A_{0}(x)=A_{0}$ the marix of $J_{0}(x)=J_{0},$ so we have $A_{0}^{2}=-I_{2},$ ($I_{2}$ is the unity matrix in dimension 2),  that is to say
\begin{equation}
A_{0}(x)=\begin{pmatrix}
a(x)&b(x)\\-\frac{1+a(x)^{2}}{b(x)}&-a(x)
\end{pmatrix}
\;\;\hbox{with}\;\; b(x)\neq 0.
\end{equation}
In the same way, we define an almost complex structure $J(y)=J$ on an open set $\widetilde{D}$ of $\mathbb{R}^{2n}.$ $J$ is $\mathbb{R}$-analytically dependind in $y\in \widetilde{D}$ and $J^{2}=-Id.$ So, the matrix $A$ associated to $J$ verifies $A^{2}=-I_{2n}$ ($I_{2n}$ is the unity matrix in dimension $2n$).\\
Then, a function 
$$f\;:\;D\longrightarrow\widetilde{D}$$
is said almost holomorphic, for the almost complex structures $(J_{0},J)$ if for all $x\in D,$ 
$$df(x)\circ J_{0}(x)=J(f(x))\circ df(x).$$
This gives, if we note $\frac{\partial f}{\partial x}(x)\in \mathcal{M}_{2n,2}(\mathbb{R})$ the matrix 
$$\frac{\partial f}{\partial x}(x)=\begin{pmatrix}
\frac{\partial f_{j}(x)}{\partial x_{i}}
\end{pmatrix}
_{j=1,...,,2n\,;\,i=1,2},$$
\begin{equation}
\frac{\partial f}{\partial x}(x)A_{0}(x)=A(f(x))\frac{\partial f}{\partial x}(x),
\end{equation}
and, therefore, the system $\forall j=1,...,2n,$
\begin{equation}
\begin{cases}
\frac{1+a^{2}}{b}\frac{\partial f_{j}}{\partial x_{2}}=(a-a_{j,j})\frac{\partial f_{j}}{\partial x_{1}}-\sum_{i\neq j}a_{j,i}\frac{\partial f_{i}}{\partial x_{1}}\\
b\frac{\partial f_{j}}{\partial x_{1}}=(a+a_{j,j})\frac{\partial f_{j}}{\partial x_{2}}+\sum_{i\neq j}a_{j,i}\frac{\partial f_{i}}{\partial x_{2}},
\end{cases}
\end{equation}
or
\begin{equation}
\begin{cases}
\frac{\partial f_{j}}{\partial x_{1}}=\frac{a+a_{j,j}}{b}\frac{\partial f_{j}}{\partial x_{2}}+\sum_{i\neq j}\frac{a_{j,i}}{b}\frac{\partial f_{i}}{\partial x_{2}}\\
\frac{\partial f_{j}}{\partial x_{2}}=\frac{b(a-a_{j,j})}{1+a^{2}}\frac{\partial f_{j}}{\partial x_{1}}-\sum_{i\neq j}\frac{ba_{j,i}}{1+a^{2}}\frac{\partial f_{i}}{\partial x_{1}}.
\end{cases}
\end{equation}
The first equation can be written, with obvious notations, 
\begin{equation}
\frac{\partial f}{\partial x_{1}}=\frac{aI+A}{b}\frac{\partial f}{\partial x_{2}},
\end{equation}
and the second one
\begin{equation}
\frac{\partial f}{\partial x_{2}}=\frac{b}{1+a^{2}}(aI-A)\frac{\partial f}{\partial x_{1}}.
\end{equation}
But, $\frac{b}{1+a^{2}}(aI-A)\frac{aI+A}{b}=I.$ Therefore, the two equations are the sames. Finally, $f$ is almost holomorphic if $\frac{\partial f}{\partial x_{2}}=\mathbb{A}\frac{\partial f}{\partial x_{1}}$ with $\mathbb{A}=\frac{b(aI-A)}{1+a^{2}}$\\
In \cite{BCGGG}, $\frac{\partial f_{j}}{\partial x_{i}}$ is noted $p^{j}_{i}$ and, therefore, the condition of almost holomorphicity is $p_{2}=\mathbb{A}p_{1}.$ \\
One of the aims of this paper is to study the possibility to have an almost holomorphic curve in a $\mathbb{R}-$analytic hypersurface. Let be $H=\{y\in \widetilde{D}:\rho(y)=0\}$ an hypersurface in $\widetilde{D}$ with $\rho$ $\mathbb{R}$-analytic. Is it possible to have $f\;:\;D\longrightarrow \widetilde{D}$ almost complex such that $\rho(D)\subset H\;?$ In other words, have we a solution $f\;:\;D\longrightarrow \widetilde{D}$ for the system of PDE 
\begin{equation}\label{B00}
\begin{cases}
\rho(f(x))=0\\
p_{2}=\mathbb{A}p_{1}\;?
\end{cases}
\end{equation}
But, in a first time, we study a more general problem. In a second time, we shall return to the almost complex case, and, also, to the complex case. \\

\section{Necessary conditions}
We are now looking for necessary conditions to have a curve $f\;:\;D\longrightarrow \widetilde{D}$  such that $\rho(D)\subset H=\{\rho=0\}$ with $\rho\;\; \mathbb{R}$-analytic and verifying $p_{2}=\mathcal{A}p_{1}$ where $\mathcal{A}=(\alpha_{j,i}(y))_{i,j=1,...,2n}\in \mathcal{M}_{2n}(\mathbb{R})$ is any matrix of order $2n$ whose the coefficients $\alpha_{j,i}$ are $\mathbb{R}$-analytic functions of $y\in D.$ Therefore, we have a solution for the system of PDE 
\begin{equation}\label{B0}
\begin{cases}
\rho(f(x))=0\\
p_{2}=\mathcal{A}p_{1}
\end{cases}
\end{equation}
which, in terms of \cite{BCGGG} (see page 131, example 5.4), is expressed by the Pfaffian differential system 
\begin{equation}\label{B1}
\begin{cases}
\rho(f(x))=0\\
p_{2}=\mathcal{A}p_{1}\\
\theta^{j}=df_{j}-p^{j}_{1}dx_{1}-p^{j}_{2}dx_{2}=0,\quad \forall j=1,...,2n, 
\end{cases}
\end{equation}
differential system which lives on the space $M_{1}$ of variables $$(x_{1}, x_{2}, f_{1},...,f_{2n},p^{j}_{i}),\;\;j=1,...,2n,\;\;i=1,2.$$ 
On account of the second line of \ref{B1}, we have $p^{j}_{2}=\sum_{i}\alpha_{j,i}p^{i}_{1},$ so we have to restrict $M_{1}$ by deleting the $p^{j}_{2}.$ And this differential system can be written
\begin{equation}\label{B2}
\begin{cases}
\rho(f(x))=0\\
\theta^{j}=df_{j}-p^{j}_{1}dx_{1}-\sum_{i}\alpha_{j,i}p^{i}_{1}dx_{2}=0,\quad \forall j=1,...,2n, 
\end{cases}
\end{equation}
on the space $M_{2}$ of variables 
$$(x_{1},x_{2},f_{1},...,f_{2n},p^{1}_{1},...,p^{2n}_{1}).$$
By deriving the first line of \ref{B2}, we obtain, if we note $\rho_{j}=\frac{\partial \rho}{\partial y_{j}}(f(x)),$ the system
\begin{equation}
\begin{cases}\label{A6}
\rho_{1}p^{1}_{1}+\rho_{2}p^{2}_{1}&=-\sum_{j=3}^{2n}\rho_{j}p^{j}_{1}\\
\sum_{j'=1}^{2n}\rho_{j'}\alpha_{j',1}p^{1}_{1}+\sum_{j'=1}^{2n}\rho_{j'}\alpha_{j',2}p^{2}_{1}&=-\sum_{j=3}^{2n}\sum_{j'=1}^{2n}\rho_{j'}\alpha_{j',j}p^{j}_{1}
\end{cases}
\end{equation}
and, when the determinant of this system $D=\sum_{j'=1}^{2n}\rho_{j'}(\rho_{1}\alpha_{j',2}-\rho_{2}\alpha_{j',1})$ is $\neq 0,$ we have 
\begin{equation}\label{A10}
\begin{cases}
p^{1}_{1}=\frac{-1}{D}\sum_{j=3}^{2n}p^{j}_{1}\sum_{j'=1}^{2n}\rho_{j'}(\rho_{j}\alpha_{j',2}-\rho_{2}\alpha_{j',j}):=\sum_{j=3}^{2n}\gamma^{1}_{j}p^{j}_{1} \\
p^{2}_{1}=\frac{-1}{D}\sum_{j=3}^{2n}p^{j}_{1}\sum_{j'=1}^{2n}\rho_{j'}(\rho_{1}\alpha_{j',j}-\rho_{j}\alpha_{j',1}):=\sum_{j=3}^{2n}\gamma^{2}_{j}p^{j}_{1}.
\end{cases}
\end{equation}
Using this formulas, the structure equations of the last line of \ref{B2} become
\begin{equation}
\begin{split}
\theta^{1}&=df_{1}-\sum_{j=3}^{2n}\gamma^{1}_{j}p^{j}_{1}dx_{1}-\sum_{j'=1}^{2n}\alpha_{1,j'}p^{j'}_{1}dx_{2}\\
&=df_{1}-\sum_{j=3}^{2n}\gamma^{1}_{j}p^{j}_{1}dx_{1}-(\alpha_{1,1}p^{1}_{1}+\alpha_{1,2}p^{2}_{1}+\sum_{j=3}^{2n}\alpha_{1,j'}p^{j}_{1})dx_{2}\\
&=df_{1}-\sum_{j=3}^{2n}\gamma^{1}_{j}p^{j}_{1}dx_{1}-\sum_{j=3}^{2n}(\alpha_{1,1}\gamma^{1}_{j}+\alpha_{1,2}\gamma^{2_{j}}+\alpha_{1,j})p^{j}_{1}dx_{2},
\end{split}
\end{equation}

\begin{equation}
\begin{split}
\theta^{2}&=df_{2}-\sum_{j=3}^{2n}\gamma^{2}_{j}p^{j}_{1}dx_{1}-\sum_{j'=1}^{2n}\alpha_{2,j'}p^{j'}_{1}dx_{2}\\
&=df_{2}-\sum_{j=3}^{2n}\gamma^{2}_{j}p^{j}_{1}dx_{1}-(\alpha_{2,1}p^{1}_{1}+\alpha_{2,2}p^{2}_{1}+\sum_{j=3}^{2n}\alpha_{2,j'}p^{j}_{1})dx_{2}\\
&=df_{2}-\sum_{j=3}^{2n}\gamma^{2}_{j}p^{j}_{1}dx_{1}-\sum_{j=3}^{2n}(\alpha_{2,1}\gamma^{1}_{j}+\alpha_{2,2}\gamma^{2_{j}}+\alpha_{2,j})p^{j}_{1}dx_{2},
\end{split}
\end{equation}
and, $\forall i=3,...,2n,$
\begin{equation}
\begin{split}
\theta^{i}&=df_{i}-p^{i}_{1}dx_{1}-\sum_{j'=1}^{2n}\alpha_{i,j'}p^{j'}_{1}dx_{2}\\
&=df_{i}-p^{i}_{1}dx_{1}-(\alpha_{_{i,1}}\sum_{j=3}^{2n}\gamma^{1}_{j}p^{j}_{1}+\alpha_{i,2}\sum_{j=3}^{2n}\gamma^{2}_{j}p^{j}_{1}+\sum_{j=3}^{2n}\alpha_{i,j}p^{j}_{1})dx_{2}\\
&=df_{i}-p^{i}_{1}dx_{1}-\sum_{j=3}^{2n}(\alpha_{i,1}\gamma^{1}_{j}+\alpha_{i,2}\gamma^{2}_{j}+\alpha_{i,j})p^{j}_{1}dx_{2}.
\end{split}
\end{equation}
We define $\beta_{i,j}=\alpha_{i,1}\gamma^{1}_{j}+\alpha_{i,2}\gamma^{2}_{j}+\alpha_{i,j}$, and, finally, the system \ref{B2} becomes
\begin{equation}
\begin{cases}\label{B4}
\theta^{1}=df_{1}-\sum_{j=3}^{2n}\gamma^{1}_{j}p^{j}_{1}dx_{1}-\sum_{j=3}^{2n}\beta_{1,j}p^{j}_{1}dx_{2}\\
\theta^{2}=df_{2}-\sum_{j=3}^{2n}\gamma^{2}_{j}p^{j}_{1}dx_{1}-\sum_{j=3}^{2n}\beta_{2,j}p^{j}_{1}dx_{2}\\
\theta^{i}=df_{i}-p^{i}_{1}dx_{1}-\sum_{j=3}^{2n}\beta_{i,j}p^{j}_{1}dx_{2}\quad \forall i=3,...,2n,
\end{cases}
\end{equation}
on the space $M$ of variables 
\begin{equation}\label{B3}
(x_{1},x_{2},f_{1},...,f_{2n},p^{3}_{1},...,p^{2n}_{1}).
\end{equation}

To solve the previous Pfaff system, we have to consider two algebraic objects, the torsion and the tableaux bundle. The idea to solve \ref{B4} is to used the following result stated in an informal way: \ref{B4} has locally solutions over a point $x$ of $M$ if the torsion vanishes locally around $x$ in $M$ and the tableau associated to \ref{B4} is in involution at $x$. The last notions will be specified in the following.\\

 From now, we calculate with the variables \ref{B3}. We remark the functions $\alpha_{i,j}, \;\gamma^{i}_{j},\;\beta_{i,j},\;\rho$ are functions of $f=(f_{1},...,f_{2n})$ only. So, $\rho_{i}=\frac{\partial \rho}{\partial y_{i}}(f)$ shall be note, sometimes, $\rho_{i}=\frac{\partial\rho}{\partial  f{i}}.$ And also, from now on, $d$ is the exterior derivative relative to the variables \ref{B3}. \\

In $T^{\ast}M,$ we take the basis 
$$\theta^{1},...,\theta^{2n},dx_{1},dx_{2},dp^{3}_{1},...,dp^{2n}_{1}.$$\\
We also consider the dual basis  noted 
$$\frac{\partial}{\partial \theta^{1}},...,\frac{\partial}{\partial \theta^{2n}},\frac{\partial}{\partial x^{1}},\frac{\partial}{\partial x^{2}},\frac{\partial}{\partial p^{3}_{1}},...,\frac{\partial}{\partial p^{2n}_{1}}$$
of $TM.$\\
\\

Let $I\subset T^{\star}M$ be the sub-bundle generated by $(\theta^{i},\; i=1,...,2n),$ $J\subset T^{\star}M$ be the sub-bundle generated by $(dx_{1}, \;dx_{2},\;\theta^{i},\; i=1,...,2n),$  and let $\{I\}\subset \Omega^{\star}(M)$ be the algebraic ideal generated by the $C^{\infty}$-sections of $I.$\\

Conformably to \cite{BCGGG}, p. 130, for $k=1,...,2n,$ we write
\begin{equation}\label{A5}
d\theta^{k}=\sum_{j=3}^{2n}\sum_{i=1}^{2}\,A^{k}_{j,i}\,dp^{j}_{1}\wedge dx_{i}+c^{k}_{1,2}\,dx_{1}\wedge dx_{2}\quad modulo \quad \{I\}.
\end{equation}
For an equality modulo $\{I\},$ we shall write $\approx.$

From \ref{B4}, we obtain 

\begin{equation}\label{A3}
\begin{cases}
df_{1}\approx\sum_{j=3}^{2n}\gamma^{1}_{j}p^{j}_{1}dx_{1}+\sum_{j=3}^{2n}\beta_{1,j}p^{j}_{1}dx_{2}\\
df_{2}\approx\sum_{j=3}^{2n}\gamma^{2}_{j}p^{j}_{1}dx_{1}+\sum_{j=3}^{2n}\beta_{2,j}p^{j}_{1}dx_{2}\\
df_{i}\approx p^{i}_{1}dx_{1}+\sum_{j=3}^{2n}\beta_{i,j}p^{j}_{1}dx_{2}\quad \forall i=3,...,2n,
\end{cases}
\end{equation}

So, using \ref{A3}, we have to calculate
\begin{equation}\label{B14}
\begin{cases}
d\theta^{1}&=-\sum_{j=3}^{2n}\sum_{i=1}^{2n}\frac{\partial}{\partial f_{i}}\Big(\gamma^{1}_{j}\Big)p^{j}_{1}df_{i}\wedge dx_{1}-\sum_{j=3}^{2n}\sum_{i=1}^{2n}\frac{\partial}{\partial f_{i}}\Big(\beta_{1,j}\Big)p^{j}_{1}df_{i}\wedge dx_{2}\\
&\quad\quad-\sum_{j=3}^{2n}\gamma^{1}_{j}dp^{j}_{1}\wedge dx_{1}-\sum_{j=3}^{2n}\beta_{1,j}dp^{j}_{1}\wedge dx_{2}\\
&\approx \sum_{j=3}^{2n}\frac{\partial \gamma^{1}_{j}}{\partial f_{1}}p^{j}_{1}\sum_{j'=3}^{2n}\beta_{1,j'}p^{j'}_{1}dx_{1}\wedge dx_{2}+\sum_{j=3}^{2n}\frac{\partial \gamma^{1}_{j}}{\partial f_{2}}p^{j}_{1}\sum_{j'=3}^{2n}\beta_{2,j'}p^{j'}_{1}dx_{1}\wedge dx_{2}\\
&\quad\quad +\sum_{j,j'=3}^{2n}\frac{\partial \gamma^{1}_{j}}{\partial f_{j'}}p^{j}_{1}\beta_{j,j'}p^{j'}_{1}dx_{1}\wedge dx_{2}-\sum_{j=3}^{2n}\gamma^{1}_{j}dp^{j}_{1}\wedge dx_{1}\\
&\qquad -\sum_{j=3}^{2n}\frac{\partial \beta_{1,j}}{\partial f_{1}}p^{j}_{1}\sum_{j'=3}^{2n}\gamma^{1}_{j'}p^{j'}_{1}dx_{1}\wedge dx_{2}-\sum_{j=3}^{2n}\frac{\partial \beta_{1,j}}{\partial f_{2}}p^{j}_{1}\sum_{j'=3}^{2n}\gamma^{2}_{j'}p^{j'}_{1}dx_{1}\wedge dx_{2}\\
&\quad\quad +\sum_{j,j'=3}^{2n}\frac{\partial \beta_{1,j}}{\partial f_{j'}}p^{j}_{1}p^{j'}_{1}dx_{1}\wedge dx_{2}-\sum_{j=3}^{2n}\beta_{1,j}dp^{j}_{1}\wedge dx_{2}\\
& \approx dx_{1}\wedge dx_{2}\sum_{j,j'=3}^{2n}p^{j}_{1}p^{j'}_{1}\Big[\frac{\partial \gamma^{1}_{j}}{\partial f_{1}}\beta_{1,j'}+\frac{\partial \gamma^{1}_{j}}{\partial f_{2}}\beta_{2,j'}+\frac{\partial \gamma^{1}_{j}}{\partial f_{j'}}\beta_{j,j'}-\frac{\partial \beta_{1,j}}{\partial f_{1}}\gamma^{1}_{j'}-\frac{\partial \beta_{1,j}}{\partial f_{2}}\gamma^{2}_{j'}-\frac{\partial \beta_{1,j}}{\partial f_{1j'}}\Big]\\
&\qquad -\sum_{j=3}^{2n}\gamma^{1}_{j}dp^{j}_{1}\wedge dx_{1}-\sum_{j=3}^{2n}\beta_{1,j}dp^{j}_{1}\wedge dx_{2}\\
d\theta^{2}&=-\sum_{j=3}^{2n}\sum_{i=1}^{2n}\frac{\partial \gamma^{2}_{j}}{\partial f_{i}}p^{j}_{1}df_{i}\wedge dx_{1}- \sum_{j=3}^{2n}\sum_{i=1}^{2n}\frac{\partial \beta_{2,j}}{\partial f_{i}}p^{j}_{1}df_{i}\wedge dx_{2}-\sum_{j=3}^{2n}\gamma^{1}_{j}dp^{j}_{1}\wedge dx_{1}\\
&\qquad -\sum_{j=3}^{2n}\beta_{2,j}dp^{j}_{1}\wedge dx_{2}\\
&\approx dx_{1}\wedge dx_{2}\sum_{j,j'=3}^{2n}p^{j}_{1}p^{j'}_{1}\Big[\frac{\partial \gamma^{2}_{j}}{\partial f_{1}}\beta_{1,j'}+\frac{\partial \gamma^{2}_{j}}{\partial f_{2}}\beta_{2,j'}+\frac{\partial \gamma^{2}_{j}}{\partial f_{j'}}\beta_{j,j'}-\frac{\partial \beta_{2,j}}{\partial f_{1}}\gamma^{1}_{j'}-\frac{\partial \beta_{2,j}}{\partial f_{2}}\gamma^{2}_{j'}-\frac{\partial \beta_{2,j}}{\partial f_{1j'}}\Big]\\
&\qquad -\sum_{j=3}^{2n}\gamma^{2}_{j}dp^{j}_{1}\wedge dx_{1}-\sum_{j=3}^{2n}\beta_{2,j}dp^{j}_{1}\wedge dx_{2}\\
&\forall i=3,...,2n,\\
d\theta^{i}&=-dp^{i}_{1}\wedge dx_{1}-\sum_{j=3}^{2n}\sum_{j'=1}^{2n}\frac{\partial \beta_{i,j}}{\partial f_{j'}}p^{j}_{1}df_{j'}\wedge dx_{2}-\sum_{j=3}^{2n}\beta_{i,j}dp^{j}_{1}\wedge dx_{2}\\
&\approx -\sum_{j,j'=3}^{2n}p^{j}_{1}p^{j'}_{1}\Big[\frac{\partial \beta_{i,j}}{\partial f_{1}}\gamma^{1}_{j'}+\frac{\partial \beta_{i,j}}{\partial f_{2}}\gamma^{2}_{j'}+\frac{\partial \beta_{i,j}}{\partial f_{j'}}\Big]dx_{1}\wedge dx_{2}-dp^{i}_{1}\wedge dx_{1}-\sum_{j=3}^{2n}\beta_{i,j}dp^{j}_{1}\wedge dx_{2}.
\end{cases}
\end{equation}
With the notation \ref{A5} adapted to our situation, we have\\
\begin{equation}\label{B5}
\begin{split}
c^{1}_{1,2}&=\sum_{j,j'=3}^{2n}p^{j}_{1}p^{j'}_{1}\Big[\frac{\partial \gamma^{1}_{j}}{\partial f_{1}}\beta_{1,j'}+\frac{\partial \gamma^{1}_{j}}{\partial f_{2}}\beta_{2,j'}+\frac{\partial \gamma^{1}_{j}}{\partial f_{j'}}\beta_{j,j'}-\frac{\partial \beta_{1,j}}{\partial f_{1}}\gamma^{1}_{j'}-\frac{\partial \beta_{1,j}}{\partial f_{2}}\gamma^{2}_{j'}-\frac{\partial \beta_{1,j}}{\partial f_{1j'}}\Big]\\ 
c^{2}_{1,2}&=\sum_{j,j'=3}^{2n}p^{j}_{1}p^{j'}_{1}\Big[\frac{\partial \gamma^{2}_{j}}{\partial f_{1}}\beta_{1,j'}+\frac{\partial \gamma^{2}_{j}}{\partial f_{2}}\beta_{2,j'}+\frac{\partial \gamma^{2}_{j}}{\partial f_{j'}}\beta_{j,j'}-\frac{\partial \beta_{2,j}}{\partial f_{1}}\gamma^{1}_{j'}-\frac{\partial \beta_{2,j}}{\partial f_{2}}\gamma^{2}_{j'}-\frac{\partial \beta_{2,j}}{\partial f_{1j'}}\Big]\\ 
c^{i}_{1,2}&=-\sum_{j,j'=3}^{2n}p^{j}_{1}p^{j'}_{1}\Big[\frac{\partial \beta_{i,j}}{\partial f_{1}}\gamma^{1}_{j'}+\frac{\partial \beta_{i,j}}{\partial f_{2}}\gamma^{2}_{j'}+\frac{\partial \beta_{i,j}}{\partial f_{j'}}\Big]\\
A^{1}_{(j,1),1}&=-\gamma^{1}_{j},\quad A^{2}_{(j,1),1}=-\gamma^{2}_{j},\quad A^{i}_{(j,1),2}=-\delta_{i}^{j},\quad  A^{1}_{(j,1),2}=-(\alpha_{1,1}\gamma^{1}_{j}+\alpha_{1,2}\gamma^{2}_{j}+\alpha_{1,j}),\\
 A^{2}_{(j,1),2}&=-(\alpha_{2,1}\gamma^{1}_{j}+\alpha_{2,2}\gamma^{2}_{j}+\alpha_{2,j}),\quad A^{i}_{(j,1),2}=-(\alpha_{i,1}\gamma^{1}_{j}+\alpha_{i,2}\gamma^{2}_{j}+\alpha_{i,j}),
 \end{split}
\end{equation}
which allow to define (see \cite{BCGGG} p. 133)
$$\pi\;:\;J^{\perp}=Span\Big(\frac{\partial}{\partial p^{3}_{1}},...,\frac{\partial}{\partial p^{2n}_{1}}\Big)\longrightarrow I^{\ast}\otimes (J/I)^{\ast}=Span\Big(\frac{\partial}{\partial \theta^{1}},...,\frac{\partial}{\partial \theta^{2n}}\Big)\otimes Span(x_{1},x_{2})$$
and we have, for 
\begin{equation}
v=\sum_{i=3}^{2n}v^{i}\frac{\partial}{\partial p^{i}_{1}}\in J^{\perp},
\end{equation}
\begin{equation}
\begin{split}
\pi(v)&=-\sum_{j=3}^{2n}\gamma^{1}_{j}v^{j}\frac{\partial}{\partial \theta^{1}}\otimes x_{1}-\sum_{j=3}^{2n}\beta_{1,j}v^{j}\frac{\partial}{\partial \theta^{1}}\otimes x_{2}\\
&\quad\quad -\sum_{j=3}^{2n}\gamma^{2}_{j}v^{j}\frac{\partial}{\partial \theta^{2}}\otimes x_{1}-\sum_{j=3}^{2n}\beta_{2,j}v^{j}\frac{\partial}{\partial \theta^{2}}\otimes x_{2}\\
&\quad\quad -\sum_{i=3}^{2n}v^{i}\frac{\partial}{\partial \theta^{i}}\otimes x_{1}-\sum_{i=3}^{2n}\sum_{j=3}^{2n}\beta_{i,j}v^{j}\frac{\partial}{\partial \theta^{i}}\otimes x_{2}\\
&=-\sum_{j=3}^{2n}v^{j}\big[\gamma^{1}_{j}\frac{\partial}{\partial \theta^{1}}\otimes x_{1}+\beta_{1,j}\frac{\partial}{\partial \theta^{1}}\otimes x_{2}+\gamma^{2}_{j}\frac{\partial}{\partial \theta^{2}}\otimes x_{1}\\
&\quad\quad +\beta_{2,j}\frac{\partial}{\partial \theta^{2}}\otimes x_{2}+\frac{\partial}{\partial\theta^{j}}\otimes x_{1}+\sum_{i=3}^{2n}\beta_{i,j}\frac{\partial}{\partial \theta^{i}}\otimes x_{2}\big]\\
&:=\sum_{j=3}^{2n}v^{j}U_{j }.
\end{split}
\end{equation}
The image of $\pi$ is the tableau $A=A^{(0)}$ associated to the system \ref{B4}, so $(U_{3},...,U_{2n})$ is a basis of $A$ (the independence of the vectors $U_{j}$ is obvious). So $dim(A)=2n-2.$\\
We also have to remark that $A_{1}=\Big\{P\in A:\frac{\partial P}{\partial x_{1}}=0\Big\}=\{0\}$ (see \cite{BCGGG} p. 119 for the definition of $A_{1}$) and therefore $dim(A_{1})=0.$\\
We now need the prolongations of $A$ (see \cite{BCGGG} p. 117). The first prolongation is 
$$A^{(1)}=\Big\{P=P_{1,1}\otimes x_{1}^{2}+P_{1,2}\otimes x_{1}x_{2}+P_{2,2}\otimes x_{2}^{2}\;:\;\frac{\partial P}{\partial x_{1}}\; and\; \frac{\partial P}{\partial x_{2}} \in A \Big\}$$
(sometimes, we shall remove the sign $\otimes$).
To obtain $P\in A^{(1)},$ we want 
\begin{equation}
\begin{cases}
\frac{\partial P}{\partial x_{1}}=2P_{1,1}x_{1}+P_{1,2}x_{2}=\sum_{j=3}^{2n}v^{j}U_{j}\\
\frac{\partial P}{\partial x_{2}}=P_{1,2}x_{1}+2P_{2,2}x_{2}=\sum_{j=3}^{2n}v'^{j}U_{j}
\end{cases}
\end{equation}
that is to say, if we explicit $U_{j},$
\begin{equation}
\begin{cases}
P_{1,1}&=\frac{1}{2}\sum_{j}v^{j}[\gamma^{1}_{j}\frac{\partial}{\partial \theta^{1}}+\gamma^{2}_{j}\frac{\partial}{\partial \theta^{2}}+\frac{\partial}{\partial \theta^{j}}]\\
P_{1,2}&=\sum_{j=3}^{2n}v^{j}[\beta_{1,j}\frac{\partial}{\partial \theta^{1}}+\beta_{2,j}\frac{\partial}{\partial \theta^{2}}+\sum_{i=3}^{2n}\beta_{i,j}\frac{\partial}{\partial \theta^{i}}]\\
&=\sum_{j=3}^{2n}v'^{j}[\gamma^{1}_{j}\frac{\partial}{\partial \theta^{1}}+\gamma^{2}_{j}\frac{\partial}{\partial \theta^{2}}+\frac{\partial}{\partial \theta^{j}}]\\
P_{2,2}&=\frac{1}{2}\sum_{j=3}^{2n}v'^{j}[\beta_{1,j})\frac{\partial}{\partial \theta^{1}}+\beta_{2,j}\frac{\partial}{\partial \theta^{2}} +\sum_{i=3}^{2n}\beta_{i,j}\frac{\partial}{\partial \theta^{i}}]
\end{cases}
\end{equation}
The equality of the second and third lines gives 
\begin{equation}\label{B6}
\begin{cases}
v'^{j}=\Sigma_{i=3}^{2n}v^{i}\beta_{j,i} \quad  \forall j=3,...,2n.\\
\sum_{j=3}^{2n}v'^{j}\gamma^{2}_{j}=\Sigma_{i=3}^{2n}v^{i}\beta_{2,i}\\
\sum_{j=3}^{2n}v'^{j}\gamma^{1}_{j}=\Sigma_{i=3}^{2n}v^{i}1,i.
\end{cases}
\end{equation}
We can translate this matricially. \\
We shall note $\alpha=\Big(\alpha_{j,i}\Big)_{j,i=3,...,2n}\in \mathcal{M}_{2n-2}(\mathbb{R}).$ $\alpha$ is the matrix obtained from $\mathcal{A}$ if we remove the two first lines and columns. \\
We also define the one column matrix $\alpha_{.,i}=\Big(\alpha_{j,i}\Big)_{j=3,...,2n}\in \mathcal{M}_{2n-2,1}(\mathbb{R})$ and the one line matrix $\gamma^{j}=\Big(\gamma^{j}_{i}\Big)_{i=3,...,2n}\in\mathcal{M}_{1,2n-2}(\mathbb{R}).$  
Let $\beta=\alpha_{.,1}\gamma^{1}+\alpha_{.,2}\gamma^{2}+\alpha =(\beta_{j,i})_{i,j=3,...,2n}\in \mathcal{M}_{2n-2}(\mathbb{R})$ be the matrix defined by 
$\beta_{j,i}=\alpha_{j,1}\gamma^{1}_{i}+\alpha_{j,2}\gamma^{2}_{i}+\alpha_{j,i},$ \\
$\beta_{2}=(\beta_{2,i})_{i=3,...,2n}\in \mathcal{M}_{1,2n-2}(\mathbb{R})$, the matrix defined by 
$\beta_{2,i}=\alpha_{2,1}\gamma^{1}_{i}+\alpha_{2,2}\gamma^{2}_{i}+\alpha_{2,i}.$ \\
and $\beta_{1}=(\beta_{1,i})_{i=3,...,2n}\in \mathcal{M}_{1,2n-2}(\mathbb{R})$, the matrix defined by 
$\beta_{1,i}=\alpha_{1,1}\gamma^{1}_{i}+\alpha_{1,2}\gamma^{2}_{i}+\alpha_{1,i}.$\\
Then, \ref{B6} can be written 
\begin{equation}
\begin{cases}
v'=\beta v\\
\gamma^{2}v'=\beta_{2} v\\
\gamma^{1}v'=\beta_{1} v,
\end{cases}
\end{equation}
which implies $(\gamma^{2}\beta-\beta_{2})v=0$ and $(\gamma^{1}\beta-\beta_{1})v=0.$\\
We also consider the linear applications $L :R^{2n-2}\longrightarrow \R^{2n-2},$  $l_{1} :\R^{ 2n-2}\longrightarrow \R,$  $l_{2} :\R^{2n-2}\longrightarrow \R,$ $\widehat{ l_{1}} :\R^{2n-2}\longrightarrow \R,$  $\widehat{ l_{2} }:\R^{2n-2}\longrightarrow \R$ respectively associated to the matrix $\beta,\;\beta_{2},\;\gamma^{2},\; \beta_{1},\;\gamma^{1},$ and \ref{B6} can be written 
\begin{equation}
\begin{cases}
v'=L( v)\\
l_{2}(v')=l_{1} (v)\\
\widehat{ l_{2}}(v')=\widehat{ l_{1}} (v),
\end{cases}
\end{equation}
and, therefore, we have $(l_{2}\circ L-l_{1})v=0$ and $(\widehat{ l_{2}}\circ L-\widehat{ l_{1}})v=0,$ that is to say, \\
$v\in Ker(l_{2}\circ L-l_{1})\cap Ker(\widehat{ l_{2}}\circ L-\widehat{ l_{1})}:=S,$ and $S$ is isomorphic to $A.$\\
It is easy to verify $A^{(1)}_{1}=\{0\}.$\\
Of course, we have $dim(A^{(1)})=dim(S)\leq 2n-2=dim(A)+dim(A_{1}).$ Following \cite{BCGGG}, p. 120, $A$ is involutive if this inequality becomes an equality, that is to say\\
\begin{equation}\label{B20}
\begin{split}
A\;\; involutive & \Leftrightarrow dim(A^{(1)})= dim(A)+dim(A_{1})\\
& \Leftrightarrow dim(A^{(1)})=dim(S)=2n-2\\
& \Leftrightarrow S=\R^{2n-2}\\
& \Leftrightarrow Ker(l_{2}\circ L-l_{1})=Ker(\widehat{ l_{2}}\circ L-\widehat{ l_{1}})=\R^{2n-2}\\
& \Leftrightarrow l_{2}\circ L=l_{1}\;\; and\;\; \widehat{ l_{2}}\circ L=\widehat{ l_{1}}\\
& \Leftrightarrow \gamma^{2}\beta=\beta_{2}\;\; and\;\; \Leftrightarrow \gamma^{1}\beta=\beta_{1}\\
& \Leftrightarrow \sum_{j=3}^{2n}\gamma^{2}_{j}\beta_{j,i}=\beta_{2,i} \quad and\quad \sum_{j=3}^{2n}\gamma^{1}_{j}\beta_{j,i}=\beta_{1,i}\qquad \forall i=3,...,2n\\
\end{split}
\end{equation}

Inductively, we define (see \cite{BCGGG}, p. 117) the q-prolongation $A^{(q)}$ of the tableau $A$ $(q\geq 1$)
\begin{equation}
A^{(q)}=\Big\{P=\sum_{l=0}^{q+1}P_{[q+1-l,l]}x_{1}^{q+1-l}x_{2}^{l}\;:\;\forall \mid J\mid =q,\;\frac{\partial P}{\partial x^{J}}\in A \Big\}.
\end{equation}
To obtain $P\in A^{(q)},$ we shall write $\frac{\partial P}{\partial x^{J}}$ in the basis $(U_{3},...,U_{2n})$ of $A.$ Let $J=[q-k,k];\;\;0\leq k\leq q$ be a multi-index with length $\mid J\mid =q$ containing $q-k$ times $1,$ and $k$ times $2.$ Then, 
\begin{equation}
\frac{\partial P}{\partial x^{J}}=\frac{\partial P}{\partial x_{1}^{q-k}\partial x_{2}^{k}}=k!(q-k+1)!P_{[q-k+1,k]}x_{1}+(q-k)!(k+1)!P_{[q-k,k+1]}x_{2}=\sum_{j=3}^{2n}v^{k}_{j}U_{j},
\end{equation}
so, $P\in A^{(q)}$ if and only if $\forall k=0,...,q,$ 
\begin{equation}
\begin{cases}
k!(q-k+1)!P_{[q-k+1,k]}&=\sum_{j=3}^{2n}v^{k}_{j}[\gamma^{1}_{j}\frac{\partial}{\partial \theta^{1}}+\gamma^{2}_{j}\frac{\partial}{\partial \theta^{2}}+\frac{\partial}{\partial \theta^{j}}]\\
(q-k)!(k+1)!P_{[q-k,k+1]}&=\sum_{j=3}^{2n}v_{j}^{k+1}\sum_{i=1}^{2n}\beta_{i,j}\frac{\partial}{\partial \theta^{i}},
\end{cases}
\end{equation}
and, therefore, 
\begin{equation}
\begin{cases}
P_{[q+1,0]}&=\frac{1}{(q+1)!}\sum_{j=3}^{2n}v^{0}_{j}[\gamma^{1}_{j}\frac{\partial}{\partial \theta^{1}}+\gamma^{2}_{j}\frac{\partial}{\partial \theta^{2}}+\frac{\partial}{\partial \theta^{j}}]\\
P_{[q-k,k+1]}&=\frac{1}{(k+1)!(q-k)!}\sum_{j=3}^{2n}v^{k+1}_{j}[\gamma^{1}_{j}\frac{\partial}{\partial \theta^{1}}+\gamma^{2}_{j}\frac{\partial}{\partial \theta^{2}}+\frac{\partial}{\partial \theta^{j}}]\\
&=\frac{1}{(k+1)!(q-k)!}\sum_{j=3}^{2n}v^{k}_{j}\sum_{i=1}^{2n}\beta_{i,j}\frac{\partial}{\partial \theta^{i}} \quad\quad \forall k=0,...,q-1\\
P_{[0,q+1]}&=\frac{1}{(q+1)!}\sum_{j=3}^{2n}v^{q}_{j}\sum_{i=3}^{2n}(\beta_{i,j}\frac{\partial}{\partial \theta^{i}}
\end{cases}
\end{equation}
So, we have, from the second and third lines, $\forall \, k=0,...,q-1,$ 
\begin{equation}
\begin{cases}
v^{k+1}=\beta v^{k}\\
\gamma^{2}v^{k+1}=\beta_{2}v^{k}\\
\gamma^{1}v^{k+1}=\beta_{1}v^{k},
\end{cases}
\end{equation}
and 
\begin{equation}
\begin{cases}
v^{k+1}=\beta^{(k+1)}v^{0}\\
\gamma^{2}\beta^{(k+1)}v^{0}=\beta_{2}\beta^{k}v^{0}\\
\gamma^{1}\beta^{(k+1)}v^{0}=\beta_{1}\beta^{k}v^{0},
\end{cases}
\end{equation}
namely
\begin{equation}\label{A26}
\begin{cases}
v^{k+1}=L^{k+1}v^{0}\\
(l_{2}\circ L-l_{1})\circ L^{k}v^{0}=0\\
(\widehat{ l_{2}}\circ L-\widehat{ l_{1}})\circ L^{k}v^{0}=0. 
\end{cases}
\end{equation}
Thus, we have $L^{k}v^{0}\in Ker(l_{2}\circ L-l_{1})\cap Ker(\widehat{ l_{2}}\circ L-\widehat{ l_{1}})=S$ or $v^{0}\in L^{-k}(S)$ and so \\ $dim(A^{(q)})=dim(\cap_{k=0}^{q-1}L^{-k}(S)).$ We verify easily $A^{(q)}_{1}=\{0\}.$\\
For each $A^{(q-1)}$ we have two possibilities : \\
either $\cap_{k=0}^{q-2}L^{-k}(S)= L^{-(q-1)}(S),$ then $dim(\cap_{k=0}^{q-1}L^{-k}(S))=dim(\cap_{k=0}^{q-2}L^{-k}(S)),$ so $dim(A^{(q)})=dim(A^{(q-1)})+dim(A^{(q-1)}_{1})$ (this later term is zero) and $A^{(q-1)}$ is involutive,\\
or $\cap_{k=0}^{q-2}L^{-k}(S)\nsubseteq L^{-(q-1)}(S),$ then $dim(\cap_{k=0}^{q-1}L^{-k}(S))<dim(\cap_{k=0}^{q-2}L^{-k}(S)),$ and $A^{(q-1)}$ is not involutive. \\
We note $\mathfrak{D_{2}}=\gamma^{2}\beta -\beta_{2}=(D_{2,3},...,D_{2,2n})$ the one line matrix of $l_{2}\circ L-l_{1},$ defined by 
\begin{equation}\label{A22}
D_{2,i}=\sum_{j=3}^{2n}\gamma^{2}_{j}\beta_{j,i}-\beta_{2,i}=\sum_{j=3}^{2n}(\alpha_{j,1}\gamma^{1}_{i}+\alpha_{j,2}\gamma^{2}_{i}+\alpha_{j,i})\gamma^{2}_{j}-(\alpha_{2,1}\gamma^{1}_{i}+\alpha_{2,2}\gamma^{2}_{i}+\alpha_{2,i}),
\end{equation}
and
$\mathfrak{D_{1}}=\gamma^{1}\beta -\beta_{1}=(D_{1,3},...,D_{1,2n})$ the one line matrix of $\widehat{l_{2}}\circ L-\widehat{ l_{1}},$ defined by 
\begin{equation}
D_{1,i}=\sum_{j=3}^{2n}\gamma^{1}_{j}\beta_{j,i}-\beta_{1,i}=\sum_{j=3}^{2n}(\alpha_{j,1}\gamma^{1}_{i}+\alpha_{j,2}\gamma^{2}_{i}+\alpha_{j,i})\gamma^{1}_{j}-(\alpha_{1,1}\gamma^{1}_{i}+\alpha_{1,2}\gamma^{2}_{i}+\alpha_{1,i}).
\end{equation}
We want to calculate $\mathfrak{D_{2}}$ more precisely. \\
We define $\mu_{i}=\sum_{j=1}^{2n}\rho_{j}\alpha_{j,i},$ and $\mu^{(2)}_{i}=\sum_{j=1}^{2n}\rho_{j}\alpha^{(2)}_{j,i},$ where $\alpha^{(2)}_{j,i}$ is the generic term of the matrix $\mathcal{A}^{2}=\big(\alpha^{(2)}_{j,i}\big)_{i,j=1,...,2n}.$ From \ref{A10}, 
\begin{equation}\label{A20}
D=\rho_{1}\mu_{2}-\rho_{2}\mu_{1},\quad \gamma^{1}_{i}=\frac{-1}{D}(\rho_{i}\mu_{2}-\rho_{2}\mu_{i})\quad and\quad  \gamma^{2}_{i}=\frac{-1}{D}(\rho_{1}\mu_{i}-\rho_{i}\mu_{1}),
\end{equation}
so we have $\forall i=3,...,2n,$
\begin{equation}
\begin{split}
\sum_{j=3}^{2n}\gamma_{j}^{2}\alpha_{j,i}&=\frac{-1}{D}\sum_{j=3}^{2n}\alpha_{j,i}\big(\rho_{1}\sum_{i'=1}^{2n}\rho_{i'}\alpha_{i',j}-\rho_{j}\sum_{i'=1}^{2n}\rho_{i'}\alpha_{i',1}\big)\\
&=\frac{-1}{D}\big[\rho_{1}\sum_{i'=1}^{2n}\rho_{i'}\sum_{j=3}^{2n}\alpha_{i',j}\alpha_{j,i}-\sum_{j=3}^{2n}\rho_{j}\alpha_{j,i}\sum_{i'=1}^{2n}\rho_{i'}\alpha_{i',1}\big]\\
&=\frac{-1}{D}\big[\rho_{1}\sum_{i'=1}^{2n}\rho_{i'}(\alpha^{(2)}_{i',i}-\alpha_{i',1}\alpha_{1,i}-\alpha_{i',2}\alpha_{2,i})-\mu_{1}(\mu_{i}-\rho_{1}\alpha_{1,i}-\rho_{2}\alpha_{2,i})\big]\\
&=\frac{-1}{D}\big[\rho_{1}\mu^{(2)}_{i}-\rho_{1}\mu_{2}\alpha_{2,i}-\mu_{1}\mu_{i}+\rho_{2}\mu_{1}\alpha_{2,i}\big].
\end{split}
\end{equation}
Therefore, 
\begin{equation}
\begin{split}
D_{2,i}&=\gamma^{1}_{i}\Big[\frac{-1}{D}\big(\rho_{1}\mu^{(2)}_{1}-\rho_{1}\mu_{2}\alpha_{2,1}-\mu_{1}^{2}+\rho_{2}\mu_{1}\alpha_{2,1}\big)\big]+\gamma^{2}_{i}\Big[\frac{-1}{D}\big(\rho_{1}\mu^{(2)}_{2}-\rho_{1}\mu_{2}\alpha_{2,2}-\mu_{1}\mu_{2}+\rho_{2}\mu_{1}\alpha_{2,2}\big)\big]\\
&\qquad\qquad -\frac{1}{D}\big(\rho_{1}\mu^{(2)}_{i}-\rho_{1}\mu_{2}\alpha_{2,i}-\mu_{1}^{i}+\rho_{2}\mu_{1}\alpha_{2,i}\big)-\alpha_{2,1}\gamma^{1}_{i}-\alpha_{2,2}\gamma^{2}_{i}-\alpha_{2,i}\\
&=\frac{-\rho_{1}}{D}\big[\mu^{(2)}_{1}\gamma^{1}_{i}+\mu^{(2)}_{2}\gamma^{2}_{i}+\mu^{(2)}_{i}\big]+\frac{\rho_{1}\mu_{2}}{D}\big[\alpha_{2,1}\gamma^{1}_{i}+\alpha_{2,2}\gamma^{2}_{i}+\alpha_{2,i}\big]+\frac{\mu_{1}}{D}\big[\mu_{1}\gamma^{1}_{i}+\mu_{2}\gamma^{2}_{i}+\mu_{i}\big]\\
&\qquad\qquad \frac{\rho_{2}\mu_{1}}{D}\big[\alpha_{2,1}\gamma^{1}_{i}+\alpha_{2,2}\gamma^{2}_{i}+\alpha_{2,i}\big]-\big[\alpha_{2,1}\gamma^{1}_{i}+\alpha_{2,2}\gamma^{2}_{i}+\alpha_{2,i}\big]\\
&=\big[\alpha_{2,1}\gamma^{1}_{i}+\alpha_{2,2}\gamma^{2}_{i}+\alpha_{2,i}\big]\big[\frac{\rho_{1}\mu_{2}-\rho_{2}\mu_{1}}{D}-1\big]\\
&\qquad\qquad+\frac{1}{D}\big\{\rho_{1}\big[\mu^{(2)}_{1}\gamma^{1}_{i}+\mu^{(2)}_{2}\gamma^{2}_{i}+\mu^{(2)}_{i}\big]+\mu_{1}\big[\mu_{1}\gamma^{1}_{i}+\mu_{2}\gamma^{2}_{i}+\mu_{i}\big]\big\}
\end{split}
\end{equation}
From \ref{A20}, in the right of the last equality, the first term is zero, so, we only have
\begin{equation}
\begin{split}
D_{2,i}&=\frac{1}{D}\big\{\rho_{1}\big[\mu^{(2)}_{1}\big(\frac{-1}{D}\big)(\rho_{i}\mu_{2}-\rho_{2}\mu_{i)}+\mu^{(2)}_{2}\big(\frac{-1}{D}\big)(\rho_{1}\mu_{i}-\rho_{i}\mu_{1)}+\mu^{(2)}_{i}\big]\\
&\qquad\qquad +\mu_{1}\big[\mu_{1}\big(\frac{-1}{D}\big)(\rho_{i}\mu_{2}-\rho_{2}\mu_{i)}+\mu_{2}\big(\frac{-1}{D}\big)(\rho_{1}\mu_{i}-\rho_{i}\mu_{1)}+\mu_{i}\big]\big\}\\
&=\frac{-1}{D^{2}}\big\{\mu^{(2)}_{1}\rho_{1}\rho_{i}\mu_{2}-\mu^{(2)}_{1}\rho_{1}\rho_{2}\mu_{i}+\mu^{(2)}_{2}\rho_{1}^{2}\mu_{i}-\mu^{(2)}_{2}\rho_{1}\rho_{i}\mu_{1}-D\rho_{1}\mu^{(2)}_{i}\\
&\qquad\qquad -\mu_{1}\mu_{i}\big[-\rho_{2}\mu_{1}+\rho_{1}\mu_{2}-D\big]\big\}.
\end{split}
\end{equation}
Always from \ref{A20}, the last term between $\{\;\}$ is zero. So, 
\begin{equation}
\begin{split}
D_{2,i}&=\frac{-\rho_{1}}{D^{2}}\big\{\mu^{(2)}_{1}\big[\rho_{i}\mu_{2}-\rho_{2}\mu_{i}\big]+\mu^{(2)}_{2}\big[\rho_{1}\mu_{i}-\rho_{i}\mu_{1}\big]+\mu^{(2)}_{i}\big[\rho_{2}\mu_{1}-\rho_{1}\mu_{2}\big]\big\}\\
&=\frac{-\rho_{1}}{D^{2}}\big\{\mu_{1}\big[\rho_{2}\mu^{(2)}_{i}-\rho_{i}\mu^{(2)}_{2}\big]+\mu_{2}\big[\rho_{i}\mu^{(2)}_{1}-\rho_{1}\mu^{(2)}_{2}\big]+\mu_{i}\big[\rho_{1}\mu^{(2)}_{2}-\rho_{2}\mu^{(2)}_{1}\big]\big\}\\
&=\frac{-\rho_{1}}{D^{2}}\big\{\rho_{1}\big[\mu_{i}\mu^{(2)}_{2}-\mu_{2}\mu^{(2)}_{i}\big]+\rho_{2}\big[\mu_{1}\mu^{(2)}_{i}-\mu_{i}\mu^{(2)}_{1}\big]+\rho_{i}\big[\mu_{2}\mu^{(2)}_{1}-\mu_{1}\mu^{(2)}_{2}\big]\big\}.
\end{split}
\end{equation}

Analogous calculations allow to express $\mathfrak{D_{1}}$ more precisely also. We obtain $\forall i=3,...,n,$
\begin{equation}
D_{1,i}=\frac{-\rho_{2}}{D^{2}}\big\{\rho_{1}\big[\mu_{i}\mu^{(2)}_{2}-\mu_{2}\mu^{(2)}_{i}\big]+\rho_{2}\big[\mu_{1}\mu^{(2)}_{i}-\mu_{i}\mu^{(2)}_{1}\big]+\rho_{i}\big[\mu_{2}\mu^{(2)}_{1}-\mu_{1}\mu^{(2)}_{2}\big]\big\}.
\end{equation}
Thus, if we note the one line matrix $\rho_{\bullet}=(\rho_{3},...,\rho_{2n}),$ $\mu_{\bullet}=(\mu_{3},...,\mu_{2n}),$ $\mu^{(2)}_{\bullet}=(\mu^{(2)}_{3},...,\mu^{(2)}_{2n}),$ we have

\begin{equation}
\begin{split}
\mathfrak{D_{1}}&=\frac{-\rho_{2}}{D^{2}}\big\{\rho_{1}\big[\mu^{(2)}_{2}\mu_{\bullet}-\mu_{2}\mu^{(2)}_{\bullet}\big]+\rho_{2}\big[\mu_{1}\mu^{(2)}_{\bullet}-\mu^{(2)}_{1}\mu_{\bullet}\big]+\big[\mu_{2}\mu^{(2)}_{1}-\mu_{1}\mu^{(2)}_{2}\big]\rho_{\bullet}\big\}:=\rho_{2}\mathfrak{D_{0}}\\
\mathfrak{D_{2}}&=\frac{-\rho_{1}}{D^{2}}\big\{\rho_{1}\big[\mu^{(2)}_{2}\mu_{\bullet}-\mu_{2}\mu^{(2)}_{\bullet}\big]+\rho_{2}\big[\mu_{1}\mu^{(2)}_{\bullet}-\mu^{(2)}_{1}\mu_{\bullet}\big]+\big[\mu_{2}\mu^{(2)}_{1}-\mu_{1}\mu^{(2)}_{2}\big]\rho_{\bullet}\big\}:=\rho_{1}\mathfrak{D_{0}} .
\end{split}
\end{equation}

So, $\rho_{1}\mathfrak{D_{1}}=\rho_{2}\mathfrak{D_{2}}$ and $\mathfrak{D_{1}}=\mathfrak{D_{2}}=0$ if and only if 
\begin{equation}\label{D2}
\rho_{1}\big[\mu_{i}\mu^{(2)}_{2}-\mu_{2}\mu^{(2)}_{i}\big]+\rho_{2}\big[\mu_{1}\mu^{(2)}_{i}-\mu_{i}\mu^{(2)}_{1}\big]+\rho_{i}\big[\mu_{2}\mu^{(2)}_{1}-\mu_{1}\mu^{(2)}_{2}\big]=-D^{2}\mathfrak{D_{0,i}}=0\qquad\quad \forall i=3,...,2n.
\end{equation}

\begin{rem}\label{C5}
Of course, $\mathfrak{D_{1}}$ and $\mathfrak{D_{2}}$ depend on  $\mathcal{A},$ so we note them $\mathfrak{D_{1}}(\mathcal{A})$ and $\mathfrak{D_{2}}(\mathcal{A}).$Using the previous expression of $\mathfrak{D_{1}}(\mathcal{A})$ and $\mathfrak{D_{2}}(\mathcal{A})),$ we prove easily, if $I$ is the unit $(2n,2n)$ matrix, $A$ an other $(2n,2n)$ matrix and $\alpha,\;\beta$ functions of $f,$ then 

\begin{equation}
\mathfrak{D}_{k}\big(\alpha(f)I+\beta(f)A\big)=\beta^{3}(f)\mathfrak{D}_{k}(A),\;\;\forall k=1,2.
\end{equation}

Moreover, if $A=(a_{i,j})_{i,j=1,...,2n}=\lambda(f)I,$ then $\mu_{A,i}=\sum_{j=1}^{2n}\rho_{j}a_{j,i}=\lambda\rho_{i}$ and $\mu^{(2)}_{A,i}=\sum_{j=1}^{2n}\rho_{j}a^{(2)}_{j,i}=\lambda^{2}\rho_{i},$ so $\mathfrak{D}_{k}(A)=0.$ \\
Therefore, if $A^{2}=\lambda(f)I,$ then $\mathcal{A}=\alpha(f)I+\beta(f)A$ verifies 
\begin{equation}\label{A24}
\mathfrak{D}_{k}(\mathcal{A})=0
\end{equation}
(because $\mu^{(2)}_{A,i}=\lambda^{2}\rho_{i}$ implies $\mathfrak{D}_{k}(\mathcal{A})=0$ from \ref{D2}).
\end{rem}

We now return after \ref{A26}. As explain there, 
$A^{(q-1)}$ is involutive if and only if \begin{equation}
\cap _{k=0}^{q-2}L^{-k}(S)= L^{-(q-1)}(S)
\end{equation}
that is (see \ref{A26})
\begin{equation}
\begin{cases}
(l_{2} \circ L-l_{1})\circ L^{k}(v)=0\quad \forall k=0,...,q-2\\
(\widehat{l} _{2} \circ L-\widehat{l} _{1})\circ L^{k}(v)=0
\end{cases}
\end{equation}
implies
\begin{equation}
\begin{cases}
 (l_{2}\circ L-l_{1})\circ L^{(q-1)}(v)=0\\
  (\widehat{l} _{2}\circ L-\widehat{l} _{1})\circ L^{(q-1)}(v)=0,
\end{cases}
\end{equation}
namely
\begin{equation}\label{A28}
\begin{cases}
(l_{2}\circ L-l_{1}) (v)=0\\
(l_{2}\circ L^{2}-l_{1}\circ L)( v)=0\\
............................... \\
................................\\
(l _{2}\circ L^{(q-1)}-l _{1}\circ L^{(q-2)})( v)=0\\
(\widehat{l} _{2}\circ L-\widehat{l} _{1}) (v)=0\\
(\widehat{l} _{2}\circ L^{2}-\widehat{l} _{1}\circ L)( v)=0\\
............................... \\
................................\\
(\widehat{l} _{2}\circ L^{(q-1)}-\widehat{l} _{1}\circ L^{(q-2)})( v)=0
\end{cases}
\end{equation}

implies
\begin{equation}\label{A29}
\begin{cases}
 (l_{2}\circ L-l_{1})\circ L^{(q-1)}(v)=0\\
  (\widehat{l} _{2}\circ L-\widehat{l} _{1})\circ L^{(q-1)}(v)=0,
\end{cases}
\end{equation}

This means that the two equations \ref{A29} are linear combinations of the  \ref{A28} ones, that is to say
\begin{equation}
\begin{cases}
(l_{2}\circ L-l_{1})\circ L^{(q-1)}=\sum_{k=0}^{q-2}\alpha_{k}(l_{2}\circ L-l_{1})\circ L^{k}+\beta_{k}(\widehat{l} _{2}\circ L-\widehat{l} _{1})\circ L^{k}\\
(\widehat{l} _{2}\circ L-\widehat{l} _{1})\circ L^{(q-1)}=\sum_{k=0}^{q-2}\alpha'_{k}(l_{2}\circ L-l_{1})\circ L^{k}+\beta'_{k}(\widehat{l} _{2}\circ L-\widehat{l} _{1})\circ L^{k}
\end{cases}
\end{equation}
or, matricially, 
\begin{equation}
\begin{cases}
\mathfrak{D_{2}}\beta^{q-1}=\alpha_{0}\mathfrak{D_{2}}+\alpha_{1}\mathfrak{D_{2}}\beta+...+\alpha_{q-2}\mathfrak{D_{2}}\beta^{q-2}+\beta_{0}\mathfrak{D_{1}}+\beta_{1}\mathfrak{D_{1}}\beta+...+\beta_{q-2}\mathfrak{D_{1}}\beta^{q-2}\\
\mathfrak{D_{1}}\beta^{q-1}=\alpha'_{0}\mathfrak{D_{2}}+\alpha'_{1}\mathfrak{D_{2}}\beta+...+\alpha'_{q-2}\mathfrak{D_{2}}\beta^{q-2}+\beta'_{0}\mathfrak{D_{1}}+\beta'_{1}\mathfrak{D_{1}}\beta+...+\beta'_{q-2}\mathfrak{D_{1}}\beta^{q-2}.
\end{cases}
\end{equation}
Finally, we obtained the following theorem and corollary:
\begin{theorem}
$A^{(q-1)}$ is involutive if 
\begin{equation}
\mathfrak{D_{0}}\beta^{q-1}\in Span(\mathfrak{D_{0}}\beta^{q-2},\;\mathfrak{D_{0}}\beta^{q-3},...,,\;\mathfrak{D_{0}}\beta,\;\mathfrak{D_{0}})
\end{equation}
\end{theorem}
\begin{cor} $A$ is involutive if (see \ref{B20})
\begin{equation}\label{B8}
\mathfrak{D_{0}}=0.
\end{equation}
\end{cor}
We can remark that, if $A^{(q_{0})}$ is involutive, then $A^{(q_{0}+k})$ is also involutive for $k\geqslant 0.$ In fact, if $\mathfrak{D_{0}}\beta^{q-1}\in Span(\mathfrak{D_{0}},\;\mathfrak{D_{0}}\beta,...,\mathfrak{D_{0}}\beta^{q-2}),$ $\mathfrak{D_{0}}\beta^{q-1}=\sum_{i=0}^{q-2}\alpha_{i}\mathfrak{D_{0}}\beta^{i},$ then $\mathfrak{D_{0}}\beta^{q}=\sum_{i=0}^{q-2}\alpha_{i}\mathfrak{D_{0}}\beta^{i+1}=\alpha_{q-2}\mathfrak{D_{0}}\beta^{q-1}+\sum_{i=0}^{q-3}\alpha_{i}\mathfrak{D_{0}}\beta^{i+1}\in Span(\mathfrak{D_{0}},\;\mathfrak{D_{0}}\beta,...,\mathfrak{D_{0}}\beta^{q-2}).$\\
Thus, we have the 
\begin{cor} $A^{(q)}$ is involutive if $q\geqslant Rank(\mathfrak{D_{0}},\;\mathfrak{D_{0}\beta},...,\mathfrak{D_{0}}\beta^{2n-3}).$ 
\end{cor}
This is the Cartan-Kuranishi theorem in the present case.\\
\\
Recall that an integral element of the system \ref{B4} is a linear subspace $E\subset T_{x}M$ such that the restriction at $E$  of all $\phi$, elements of the differential ideal $I$, vanishes on $E$. The torsion is a necessary and sufficient condition such that \ref{B4} has an integral element over $x\in M$ and we calculate it below.

With the notations of \cite{BCGGG}, we define (see \ref{B5}) 
\begin{equation}
c=\Big(c^{1}_{1,2}\frac{\partial}{\partial \theta^{1}}+c^{2}_{1,2}\frac{\partial}{\partial \theta^{2}}+\sum_{i=3}^{2n}c^{i}_{1,2}\frac{\partial}{\partial \theta^{i}}\Big)\otimes dx_{1}\wedge dx_{2}
\end{equation}
which is a section of the bundle $I^{\ast}\otimes \bigwedge^{2}(J/I)$ and has to be quotiented by $Image(\overline{\pi})$ where 
\begin{equation}
\begin{split}
&\overline{\pi}\;:\;J^{\perp}\otimes J/I=Span\Big(\frac{\partial}{\partial p^{3}_{1}},...,\frac{\partial}{\partial p^{2n}_{1}}\Big)\otimes Span(dx_{1},dx_{2})\\
&\qquad\qquad\longrightarrow I^{\ast}\wedge^{2}(J/I)=Span\Big(\frac{\partial}{\partial \theta^{1}},...,\frac{\partial}{\partial \theta^{2n}}\Big)\otimes Span(dx_{1}\wedge dx_{2})
\end{split}
\end{equation}
is defined by
\begin{equation}
v=\sum_{i=3}^{2n}\sum_{k=1}^{2}v^{i}_{k}\frac{\partial}{\partial p^{i}_{1}}\otimes dx_{k}\longrightarrow \overline{\pi}(v)=\sum_{j=2}^{2n}\big(\sum_{i=3}^{2n}A^{j}_{(i,1),1}v^{i}_{2}-A^{j}_{(i,1),2}v^{i}_{1}\big)\frac{\partial}{\partial \theta^{j}}\otimes dx_{1}\wedge dx_{2}
\end{equation}
and therefore, by \ref{B5},
\begin{equation}\begin{split}
\overline{\pi}(v)&=\Big[\Big(\sum_{i=3}^{2n}\gamma^{1}_{i}v^{i}_{2}-\beta_{1,i}v^{i}_{1}\Big)\frac{\partial}{\partial \theta^{1}}+ \Big(\sum_{i=3}^{2n}\gamma^{2}_{i}v^{i}_{2}-\beta_{2,i}v^{i}_{1}\Big)\frac{\partial}{\partial \theta^{2}}\\
&\qquad\qquad +\sum_{j=3}^{2n}\Big(v^{j}_{2}-\sum_{i=3}^{2n}\beta_{j,i}v^{i}_{1}\Big)\frac{\partial}{\partial \theta^{j}}\Big]\otimes dx_{1}\wedge dx_{2}.
\end{split}
\end{equation}
Thus, the torsion vanishes if there exists $v$ such that $\overline{\pi}(v)=c,$ i.e.
\begin{equation}
\begin{cases}\label{A14}
\sum_{i=3}^{2n}\gamma^{1}_{i}v^{i}_{2}-\beta_{1,i}v^{i}_{1}=c^{1}_{1,2}\\
\sum_{i=3}^{2n}\gamma^{2}_{i}v^{i}_{2}-\beta_{2,i}v^{i}_{1}=c^{2}_{1,2}\\
v^{j}_{2}-\sum_{i=3}^{2n}\beta_{j,i}v^{i}_{1}=c^{j}_{1,2} \qquad \forall j=3,...,2n.
\end{cases}
\end{equation}
From the last lines, we obtain, $\forall j=3,...,2n,$ 
\begin{equation}
v^{j}_{2}=c^{j}_{1,2}+\sum_{i=3}^{2n}\beta_{j,i}v^{i}_{1},
\end{equation}
and we carry over the first lines
\begin{equation}
\begin{cases}
\sum_{i=3}^{2n}\gamma^{1}_{i}[c^{i}_{1,2}+\sum_{i'=3}^{2n}\beta_{i,i'}v^{i'}_{1}]-\beta_{1,i}v^{i}_{1}=c^{1}_{1,2}\\
\sum_{i=3}^{2n}\gamma^{2}_{i}[c^{i}_{1,2}+\sum_{i'=3}^{2n}\beta_{i,i'}v^{i'}_{1}]-\beta_{2,i}v^{i}_{1}=c^{2}_{1,2},
\end{cases}
\end{equation}
or
\begin{equation}
\begin{cases}
\sum_{i'=3}^{2n}v^{i'}_{1}\big[\sum_{i=3}^{2n}\beta_{i,i'}\gamma^{1}_{i}-\beta_{1,i'}\big]=c^{1}_{1,2}-\sum_{i=3}^{2n}\gamma^{1}_{i}c^{i}_{1,2}\\
\sum_{i'=3}^{2n}v^{i'}_{1}\big[\sum_{i=3}^{2n}\beta_{i,i'}\gamma^{2}_{i}-\beta_{2,i'}\big]=c^{2}_{1,2}-\sum_{i=3}^{2n}\gamma^{2}_{i}c^{i}_{1,2},
\end{cases}
\end{equation}
that is, matricialy, 
\begin{equation}\label{B9}
\begin{cases}
\mathfrak{D_{1}}v_{1}=c^{1}_{1,2}-\sum_{i=3}^{2n}\gamma^{1}_{i}c^{i}_{1,2}\\
\mathfrak{D_{2}}v_{1}=c^{2}_{1,2}-\sum_{i=3}^{2n}\gamma^{2}_{i}c^{i}_{1,2},
\end{cases}
\end{equation}

that is to say 
\begin{equation}
\begin{cases}
\rho_{2}\mathfrak{D_{0}}v_{1}=c^{1}_{1,2}-\sum_{i=3}^{2n}\gamma^{1}_{i}c^{i}_{1,2}\\
\rho_{1}\mathfrak{D_{0}}v_{1}=c^{2}_{1,2}-\sum_{i=3}^{2n}\gamma^{2}_{i}c^{i}_{1,2},
\end{cases}
\end{equation}
which implies 
\begin{equation}\label{A30}
\rho_{1}\Big[c^{1}_{1,2}-\sum_{i=3}^{2n}\gamma^{1}_{i}c^{i}_{1,2}\Big]=\rho_{2}\Big[c^{2}_{1,2}-\sum_{i=3}^{2n}\gamma^{2}_{i}c^{i}_{1,2}\Big]
\end{equation}
To summarize, we have the theorem:
\begin{theorem}\label{C6}
We have two possibilities : \\
- if $\mathfrak{D_{0}}\neq 0$ the torsion can be absorbed if \ref{A30} is satisfied, \\
-if $\mathfrak{D_{0}}=0$ if
\begin{equation}\label{A32}
\begin{cases}
c^{1}_{1,2}=\sum_{i=3}^{2n}\gamma^{1}_{i}c^{i}_{1,2}\\
c^{2}_{1,2}=\sum_{i=3}^{2n}\gamma^{2}_{i}c^{i}_{1,2}.
\end{cases}
\end{equation}
\end{theorem}
In the almost complex case, the previous theorem gives immediately 
\begin{cor}  The condition $\mathfrak{D_{0}}=0$ is satisfied in the almost complex case. So, in almost complex analysis, $A^{(0)}$ is involutive, and the torsion may be absorbed if \ref{A32} is verified. 
\end{cor}
In fact, in the almost complex case, as explain in the paragraph 2, we have $\mathcal{A}=\mathbb{A}=\frac{b(aI-A)}{1+a^{2}},$ and the remark \ref{C5} gives the corollary. \qed
\medskip

In the complex case, we can specify a little more 
\begin{cor} In complex analysis, we have $\mathfrak{D_{0}}=0$ and the torsion may be absorbed if \\
$c^{1}_{1,2}=c^{2}_{1,2}=0.$
\end{cor}
In the complex case, we have $f\;:\;D\longrightarrow \widetilde{D}$ holomorphic i.e. \\$\forall i=1,...,n,\quad p^{2i-1}_{2}=-p^{2i}_{1}\;and\; p^{2i}_{2}=p^{2i-1}_{1}$ with $\rho(f(x))=0.$\
Therefore, we have, \\
for $i=1,...,n,$
\begin{equation}
\begin{cases}
\alpha_{2i-1,2i}=-1\\
\alpha_{2i,2i-1}=1\\
\alpha_{j,j'}=0\; in\;the \;other\; cases.
\end{cases}
\end{equation}

Deriving $\rho(f(x))=0$ relatively to $x_{1}$ and $x_{2}$ gives
\begin{equation}
\begin{cases}
\rho_{1}p^{1}_{1}+\rho_{2}p^{2}_{1}=-\sum_{j=2}^{n}(\rho_{2j-1}p^{2j-1}_{1}+\rho_{2j}p^{2j}_{1})\\
\rho_{2}p^{1}_{1}-\rho_{1}p^{2}_{1}=-\sum_{j=2}^{n}(\rho_{2j}p^{2j-1}_{1}-\rho_{2j-1}p^{2j}_{1}),
\end{cases}
\end{equation}
system whose the determinant $D=-(\rho^{2}_{1}+\rho^{2}_{2})$ is supposed non zero, from which
\begin{equation}\label{B22}
\begin{cases}
p^{1}_{1}=\Sigma_{j=2}^{n}\Big(p^{2j-1}_{1}\frac{\rho_{2}\rho_{2j}+\rho_{1}\rho_{2j-1}}{D}+p^{2j}_{1}\frac{\rho_{1}\rho_{2j}-\rho_{2}\rho_{2j-1}}{D}\Big)\\
p^{2}_{1}=\Sigma_{j=2}^{n}\Big(p^{2j-1}_{1}\frac{\rho_{2}\rho_{2j-1}-\rho_{1}\rho_{2j}}{D}+p^{2j}_{1}\frac{\rho_{1}\rho_{2j-1}+\rho_{2}\rho_{2j}}{D}\Big).
\end{cases}
\end{equation}
So, we have
\begin{equation}
\begin{cases}
\gamma^{1}_{2j-1}=\gamma^{2}_{2j}=\frac{\rho_{1}\rho_{2j-1}+\rho_{2}\rho_{2j}}{D}:=\frac{D_{2,j}}{D}\\
\gamma^{1}_{2j}=-\gamma^{2}_{2j-1}=\frac{\rho_{1}\rho_{2j}-\rho_{2}\rho_{2j-1}}{D}:=\frac{D_{1,j}}{D}\\
\beta_{1,j}=-\gamma^{2}_{j}\quad and \quad \beta_{2,j}=\gamma^{1}_{j}\\
\beta_{2i-1,j}=-\delta_{2i}^{j}\quad and \quad \beta_{2i,j}=\delta_{2i-1}^{j}.
\end{cases}
\end{equation}

From \ref{B4} and \ref{B22}, we obtain 
\begin{equation}\label{C2}
\begin{cases}
\theta^{1}&=df_{1}-\sum_{j=3}^{2n}\gamma^{1}_{j}p^{j}_{1}dx_{1}+\sum_{j=3}^{2n}\gamma^{2}_{j}p^{j}_{1}dx_{2}\\
\theta^{2}&=df_{2}-\sum_{j=3}^{2n}\gamma^{2}_{j}p^{j}_{1}dx_{1}-\sum_{j=3}^{2n}\gamma^{1}_{j}p^{j}_{1}dx_{2}\\
&\forall j=2,...,n\\
\theta_{2j-1}&=df_{2j-1}-p^{2j-1}_{1}dx_{1}+p^{2j}_{1}dx_{2}\\
\theta_{2j}&=df_{2j}-p^{2j}_{1}dx_{1}-p^{2j-1}_{1}dx_{2}.
\end{cases}
\end{equation}
Consequently, 
\begin{equation}
\forall j=2,...,n,\quad d\theta^{2j-1}=-dp^{2j-1}_{1}\wedge dx_{1}+dp^{2j}_{1}\wedge dx_{2}\quad and\quad d\theta^{2j}=-dp^{2j}_{1}\wedge dx_{1}-dp^{2j-1}_{1}\wedge dx_{2},
\end{equation}
and so, $c^{2j}_{1,2}=c^{2j-1}_{1,2}=0\quad \forall j=2,...,n,$ and, therefore the torsion may be absorbed if  $c^{1}_{1,2}=c^{2}_{1,2}=0.$ \qed
\begin{rem}
The two last conditions in almost and complex analysis are  nothing else than the usual Levi form in an intrinsic way. In this situation the torsion (or Levi form) contains all the compatibility conditions to solve the system \ref{B4} because the tableau associated is in involution. It  is the unique  reason, in complex analysis or almost complex analysis, that the Levi form plays a crucial role. This is not the case when $\mathfrak{D_{0}}\not =0$ (see the paragraph successive torsion).
\end{rem}
After this remark, we, now, continue to study the complex case. \\
From the two first lines of \ref{C2}, we have
\begin{equation}\label{A31}
\begin{cases}
d\theta^{1}&=-\sum_{j=3}^{2n}\gamma^{1}_{j}dp^{j}_{1}\wedge dx_{1}+\sum_{j=3}^{2n}\gamma^{2}_{j}dp^{j}_{1}\wedge dx_{2}\\
&\qquad -\sum_{j=3}^{2n}\sum_{i=1}^{2n}\Big(p^{j}_{1}\frac{\partial \gamma^{1}_{j}}{\partial f_{i}}df_{i}\wedge dx_{1}-p^{j}_{1}\frac{\partial \gamma^{2}_{j}}{\partial f_{i}}df_{i}\wedge dx_{2}\Big)\\
&\approx -\sum_{j=3}^{2n}\gamma^{1}_{j}dp^{j}_{1}\wedge dx_{1}+\sum_{j=3}^{2n}\gamma^{2}_{j}dp^{j}_{1}\wedge dx_{2}\\
&\qquad -dx_{1}\wedge dx_{2}\Big[\sum_{j=3}^{2n}p^{j}_{1}\Big(\frac{\partial \gamma^{1}_{j}}{\partial f_{1}}\sum_{j'=3}^{2n}\gamma^{2}_{j'}p^{j'}_{1}-\frac{\partial \gamma^{1}_{j}}{\partial f_{2}}\sum_{j'=3}^{2n}\gamma^{1}_{j'}p^{j'}_{1}\\
&\qquad  +\sum_{i=2}^{n}\big(\frac{\partial \gamma^{1}_{j}}{\partial f_{2i-1}}p^{2i}_{1}-\frac{\partial \gamma^{1}_{j}}{\partial f_{2i}}p^{3i-1}_{1}\big) -\frac{\partial \gamma^{2}_{j}}{\partial f_{1}}\sum_{j'=3}^{2n}\gamma^{1}_{j'}p^{j'}_{1}-\frac{\partial \gamma^{2}_{j}}{\partial f_{2}}\sum_{j'=3}^{2n}\gamma^{2}_{j'}p^{j'}_{1}\\
&\qquad -\sum_{j'=3}^{2n}\frac{\partial \gamma^{2}_{j}}{\partial f_{j'}}p^{j'}_{1}\Big)\Big]\\
d\theta^{2}&=-\sum_{j=3}^{2n}\gamma^{2}_{j}dp^{j}_{1}\wedge dx_{1}-\sum_{j=3}^{2n}\gamma^{1}_{j}dp^{j}_{1}\wedge dx_{2}\\
&\qquad -\sum_{j=3}^{2n}\sum_{i=1}^{2n}\Big(p^{j}_{1}\frac{\partial \gamma^{2}_{j}}{\partial f_{i}}df_{i}\wedge dx_{1}+p^{j}_{1}\frac{\partial \gamma^{1}_{j}}{\partial f_{i}}df_{i}\wedge dx_{2}\Big)\\
&\approx -\sum_{j=3}^{2n}\gamma^{2}_{j}dp^{j}_{1}\wedge dx_{1}-\sum_{j=3}^{2n}\gamma^{1}_{j}dp^{j}_{1}\wedge dx_{2}\\
&\qquad -dx_{1}\wedge dx_{2}\Big[\sum_{j=3}^{2n}p^{j}_{1}\Big(\frac{\partial \gamma^{2}_{j}}{\partial f_{1}}\sum_{j'=3}^{2n}\gamma^{2}_{j'}p^{j'}_{1}-\frac{\partial \gamma^{2}_{j}}{\partial f_{2}}\sum_{j'=3}^{2n}\gamma^{1}_{j'}p^{j'}_{1}\\
&\qquad +\sum_{i=2}^{n}\big(\frac{\partial \gamma^{2}_{j}}{\partial f_{2i-1}}p^{2i}_{1}-\frac{\partial \gamma^{2}_{j}}{\partial f_{2i}}p^{3i-1}_{1}\big)+\frac{\partial \gamma^{1}_{j}}{\partial f_{1}}\sum_{j'=3}^{2n}\gamma^{1}_{j'}p^{j'}_{1}+\frac{\partial \gamma^{1}_{j}}{\partial f_{2}}\sum_{j'=3}^{2n}\gamma^{2}_{j'}p^{j'}_{1}\\
&\qquad +\sum_{j'=3}^{2n}\frac{\partial \gamma^{1}_{j}}{\partial f_{j'}}p^{j'}_{1}\Big)\Big]\\
\end{cases}
\end{equation}
and, by a long but easy calculation, we obtain
\begin{equation}
\begin{split}
c^{1}_{1,2}&=\sum_{j,j'=2}^{n}[p^{2j-1}_{1}p^{2j'-1}_{1}+p^{2j}_{1}p^{2j'}_{1}][P^{2}_{j'}(\gamma^{2}_{2j})-P^{1}_{j'}(\gamma^{1}_{2j})]\\
&\qquad +[p^{2j}_{1}p^{2j'-1}_{1}-p^{2j-1}_{1}p^{2j'}_{1}][P^{2}_{j'}(\gamma^{1}_{2j})+P^{1}_{j'}(\gamma^{2}_{2j})]
\end{split}
\end{equation}
and
\begin{equation}
\begin{split}
c^{2}_{1,2}&=\sum_{j,j'=2}^{n}-[p^{2j-1}_{1}p^{2j'-1}_{1}+p^{2j}_{1}p^{2j'}_{1}][P^{2}_{j'}(\gamma^{1}_{2j})+P^{1}_{j'}(\gamma^{2}_{2j})]\\
&\qquad +[p^{2j}_{1}p^{2j'-1}_{1}-p^{2j-1}_{1}p^{2j'}_{1}][P^{2}_{j'}(\gamma^{2}_{2j})-P^{1}_{j'}(\gamma^{1}_{2j})]
\end{split}
\end{equation}
where, for all $k=2,...,n,$ we have introduced the differential polynomials 
\begin{equation}
\begin{cases}
P^{1}_{k}=\gamma^{2}_{2k}\frac{\partial }{\partial y_{1}}-\gamma^{1}_{2k}\frac{\partial }{\partial y_{2}}+\frac{\partial }{\partial y_{2k-1}}\\
P^{2}_{k}=\gamma^{1}_{2k}\frac{\partial }{\partial y_{1}}+\gamma^{2}_{2k}\frac{\partial }{\partial y_{2}}+\frac{\partial }{\partial y_{2k}}.
\end{cases}
\end{equation}

Thus, if we note, for $j,k=2,...,n,$ 
\begin{equation}
\begin{split}
B_{j,k}&=P^{2}_{k}(\gamma^{1}_{2j})+P^{1}_{k}(\gamma^{2}_{2j)}\\
B^{j,k}&=P^{2}_{k}(\gamma^{2}_{2j})-P^{1}_{k}(\gamma^{1}_{2j}),
\end{split}
\end{equation}

then 
\begin{equation}\label{A11}
\begin{cases}
c^{1}_{1,2}&=\sum_{j,j'=2}^{2n[}[p^{2j-1}_{1}p^{2j'-1}_{1}+p^{2j}_{1}p^{2j'}_{1}]B^{j,j'}+[p^{2j}_{1}p^{2j'-1}_{1}-p^{2j-1}_{1}p^{2j'}_{1}]B_{j,j'}\\
c^{2}_{1,2}&=\sum_{j,j'=2}^{2n[}-[p^{2j-1}_{1}p^{2j'-1}_{1}+p^{2j}_{1}p^{2j'}_{1}]B_{j,j'}+[p^{2j}_{1}p^{2j'-1}_{1}-p^{2j-1}_{1}p^{2j'}_{1}]B^{j,j'}.
\end{cases}
\end{equation}

We have, here, two quadratic polynomials in $p^{j}_{1}$ which have to be $0$ in order that the torsion vanishes. If   one of the two polynomials is non degenerated positive or non degenerated negative, then the only solution is $p^{j}_{1}=0.$ Therefore, the only holomorphic discs in $H$ are points. We have no proper holomorphic disc in $H.$\\
To have really an holomorphic disc in $H,$ we need, for our quadratic polynomials, non zero solutions.

\begin{ex}
The complex case in dimension 6.
\end{ex}
We have, with the previous notations, 
\begin{equation}
\begin{split}
c^{1}_{1,2}&=B^{2,2}\big[(p^{3}_{1})^{2}+(p^{4}_{1})^{2}\big]+B^{3,3}\big[(p^{5}_{1})^{2}+(p^{6}_{1})^{2}\big]+\big[p^{3}_{1}p^{5}_{1}+p^{4}_{1}p^{6}_{1}\big]\big(B^{2,3}+B^{3,2}\big)\\
&\qquad\qquad +\big[p^{4}_{1}p^{5}_{1}-p^{3}_{1}p^{6}_{1}\big]\big(B_{2,3}-B_{3,2}\big)\\
&=B^{2,2}\Big[\Big(p^{3}_{1}+\frac{B^{2,3}+B^{3,2}}{2B_{2,2}}p^{5}_{1}-\frac{B_{2,3}-B_{3,2}}{2B_{2,2}}p^{6}_{1}\Big)^{2}+
\Big(p^{4}_{1}+\frac{B^{2,3}+B^{3,2}}{2B_{2,2}}p^{6}_{1}+\frac{B_{2,3}-B_{3,2}}{2B_{2,2}}p^{5}_{1}\Big)^{2}\Big]\\
&\qquad\qquad+\frac{4B_{2,2}B_{3,3}-(B_{2,3}-B_{3,2})^{2}-(B^{2,3}+B^{3,2})^{2}}{4B_{2,2}}\Big[(p^{5}_{1})^{2}+(p^{6}_{1})^{2}\Big]\\
c^{2}_{1,2}&=-B_{2,2}\big[(p^{3}_{1})^{2}+(p^{4}_{1})^{2}\big]-B_{3,3}\big[(p^{5}_{1})^{2}+(p^{6}_{1})^{2}\big]-\big[p^{3}_{1}p^{5}_{1}+p^{4}_{1}p^{6}_{1}\big]\big(B_{2,3}+B_{3,2}\big)\\
&\qquad\qquad +\big[p^{4}_{1}p^{5}_{1}-p^{3}_{1}p^{6}_{1}\big]\big(B^{2,3}-B^{3,2}\big)\\
&=-B_{2,2}\Big[\Big(p^{3}_{1}+\frac{B_{2,3}+B_{3,2}}{2B_{2,2}}p^{5}_{1}+\frac{B^{2,3}-B^{3,2}}{2B_{2,2}}p^{6}_{1}\Big)^{2}+
\Big(p^{4}_{1}+\frac{B_{2,3}+B_{3,2}}{2B_{2,2}}p^{6}_{1}-\frac{B^{2,3}-B^{3,2}}{2B_{2,2}}p^{5}_{1}\Big)^{2}\Big]\\
&\qquad\qquad+\frac{-4B_{2,2}B_{3,3}+(B_{2,3}+B_{3,2})^{2}+(B^{2,3}-B^{3,2})^{2}}{4B_{2,2}}\Big[(p^{5}_{1})^{2}+(p^{6}_{1})^{2}\Big].
\end{split}
\end{equation}
Now, the existence of a non trivial (non constant) solution for \ref{B0} implies that the torsion $c^{1}_{1,2}=c^{2}_{1,2}=0$ with some $p^{i}_{j}\neq 0.$ So, the quadratic forms $c^{2}_{1,2}$ and$c^{2}_{1,2},$  quadratic forms in $p^{i}_{j},$ are not definite positive (or negative). Therefore 
\begin{equation}\label{B10}
\begin{cases}
4B_{2,2}B_{3,3}-(B_{2,3}+B_{3,2})^{2}-(B^{2,3}-B^{3,2})^{2}\leq 0\\
4B^{2,2}B^{3,3}-(B^{2,3}+B^{3,2})^{2}-(B_{2,3}-B_{3,2})^{2}\leq 0.
\end{cases}
\end{equation}
An other method would be to look for an orthogonal basis for the quadratic form $c^{2}_{1,2}.$ This, of course, gives the same result. 
\begin{ex}
The pseudo-ellipso\"\i ds
\end{ex}
To end with the complex case, we give an example : suppose now we have for $\rho$ the polynomial 
\begin{equation}
\rho(y)=\alpha_{1}y_{1}^{2k_{1}}+\alpha_{2}y_{2}^{2k_{2}}+\alpha_{3}y_{3}^{2k_{3}}+\alpha_{4}y_{4}^{2k_{4}}+\alpha_{5}y_{5}^{2k_{5}}+\alpha_{6}y_{6}^{2k_{6}},
\end{equation}
with $\alpha_{1},...,\alpha_{6}\in \R,$ and $k_{1},...,k_{6}\in \N.$ \\
We note  $v_{i}=\rho_{i}=2\alpha_{i}k_{i}y_{i}^{2k_{i}-1},\quad i=1,...,6,$ and also $w_{i}=\frac{\partial v_{i}}{\partial y_{i}}=\frac{\partial^{2}\rho}{\partial^{2}y_{i}}=2k_{i}(2k_{i}-1)\alpha_{i}y_{i}^{2k_{i}-2},$ and we want to explicit the coefficients $B_{j,k}$ and $B^{j,k}$ of \ref{B10}. \\

We have, with the previous definitions applied to this example, $D=-(v_{1}^{2}+v_{2}^{2}),\;\;\gamma^{1}_{2j-1}=\gamma^{2}_{2j}=\frac{\rho_{1}\rho_{2j-1}+\rho_{2}\rho_{2j}}{D}=\frac{D_{2,j}}{D},\;\; \gamma^{1}_{2j}=-\gamma^{2}_{2j-1}=\frac{\rho_{1}\rho_{2j}-\rho_{2}\rho_{2j-1}}{D}=\frac{D_{1,j}}{D}$ and 
\begin{equation}
\begin{split}
B_{2,3}&=P^{2}_{3}(\gamma^{1}_{4})+P^{1}_{3}(\gamma^{2}_{4})\\
&=\frac{1}{D^{2}}\Big[\gamma^{1}_{6}\Big(\frac{\partial D_{1,2}}{\partial y_{1}}D-D_{1,2}\frac{\partial D}{\partial y_{1}}\Big)+\gamma^{2}_{6}\Big(\frac{\partial D_{1,2}}{\partial y_{2}}D-D_{1,2}\frac{\partial D}{\partial y_{2}}\Big)+\frac{\partial D_{1,2}}{\partial y_{6}}D-D_{1,2}\frac{\partial D}{\partial y_{6}}\\
&\qquad +\gamma^{2}_{6}\Big(\frac{\partial D_{2,2}}{\partial y_{1}}D-D_{2,2}\frac{\partial D}{\partial y_{1}}\Big)-\gamma^{1}_{6}\Big(\frac{\partial D_{2,2}}{\partial y_{2}}D-D_{2,2}\frac{\partial D}{\partial y_{2}}\Big)+\frac{\partial D_{2,2}}{\partial y_{5}}D-D_{1,2}\frac{\partial D}{\partial y_{5}}\Big]\\
&=\frac{1}{D^{3}}\{(v_{1}v_{6}-v_{2}v_{5})[-v_{4}w_{1}(v_{1}^{2}+v_{2}^{2})+(v_{1}v_{4}-v_{2}v_{3})2v_{1}w_{1}]\\
&\qquad +(v_{1}v_{5}+v_{2}v_{6})[v_{3}w_{1}(v_{2}^{2}+v_{2}^{2})+(v_{1}v_{4}-v_{2}v_{3})2v_{2}w_{2}]\\
&\qquad +(v_{1}v_{5}+v_{2}v_{6})[-v_{3}w_{1}(v_{1}^{2}+v_{1}^{2})+(v_{1}v_{3}+v_{2}v_{4})2v_{1}w_{1}]\\
&\qquad -(v_{1}v_{6}-v_{2}v_{5})[-v_{4}w_{2}(v_{1}^{2}+v_{1}^{2})+(v_{1}v_{3}+v_{2}v_{4})2v_{2}w_{2}]\}\\
&=\frac{1}{D^{3}}\{w_{1}[(v_{1}^{2}+v_{2}^{2})\big(-v_{4}(v_{1}v_{6}-v_{2}v_{5})-v_{3}(v_{1}v_{5}+v_{2}v_{6})\big)\\
&\qquad \qquad +2v_{1}\big((v_{1}v_{6}-v_{2}v_{5})(v_{1}v_{4}-v_{2}v_{3})+(v_{1}v_{5}+v_{2}v_{6})(v_{1}v_{3}+v_{2}v_{4})\big)]\\
&\qquad\quad +w_{2}[(v_{1}^{2}+v_{2}^{2})\big(v_{3}(v_{1}v_{5}+v_{2}v_{6})+v_{4}(v_{1}v_{6}-v_{2}v_{5})\big)\\
&\qquad \qquad +2v_{2}\big((v_{1}v_{5}+v_{2}v_{6})(v_{1}v_{4}-v_{2}v_{3})+(v_{1}v_{6}-v_{2}v_{5})(v_{1}v_{3}+v_{2}v_{4})\big)]\\
&=\frac{1}{D^{2}}[v_{2}(v_{3}v_{6}-v_{4}v_{5})-v_{1}(v_{3}v_{5}+v_{4}v_{6})][w_{1}+w_{2}].
\end{split}
\end{equation}

Analogous calculations give
\begin{equation}\label{A15}
\begin{split}
B_{3,2}&=\frac{1}{D^{2}}[-v_{1}(v_{3}v_{5}+v_{4}v_{6})-v_{2}(v_{3}v_{6}-v_{4}v_{5})][w_{1}+w_{2}]\\
B^{2,3}&=\frac{1}{D^{2}}[-v_{1}(v_{3}v_{6}-v_{4}v_{5})-v_{2}(v_{3}v_{5}+v_{4}v_{6})][w_{1}+w_{2}]\\
B^{3,2}&=\frac{1}{D^{2}}[v_{1}(v_{3}v_{6}-v_{4}v_{5})-v_{2}(v_{3}v_{5}+v_{4}v_{6})][w_{1}+w_{2}]\\
B_{2,2}&=\frac{-v_{1}}{D^{2}}[(v_{3}^{2}+v_{4}^{2}) ( w_{1}+w_{2})+ (v_{1}^{2}+v_{2}^{2}) ( w_{3}+w_{4}) ] \\
B_{3,3}&=\frac{-v_{1}}{D^{2}}[(v_{5}^{2}+v_{6}^{2}) ( w_{1}+w_{2})+ (v_{1}^{2}+v_{2}^{2}) ( w_{5}+w_{6}) ] \\
B^{2,2}&=\frac{-v_{2}}{D^{2}}[(v_{3}^{2}+v_{4}^{2}) ( w_{1}+w_{2})+ (v_{1}^{2}+v_{2}^{2}) ( w_{3}+w_{4}) ] \\
B^{3,3}&=\frac{-v_{2}}{D^{2}}[(v_{5}^{2}+v_{6}^{2}) ( w_{1}+w_{2})+ (v_{1}^{2}+v_{2}^{2}) ( w_{5}+w_{6}) ] \\
\end{split}
\end{equation}
So, the first equation of \ref{B10} becomes
\begin{equation}
\begin{split}
&\frac{4v_{1}^{2}(v_{1}^{2}+v_{2}^{2})}{D^{4}}[(v_{3}^{2}+v_{4}^{2})(w_{1}+w_{2})(w_{5}+w_{6)}+(v_{5}^{2}+v_{6}^{2})(w_{1}+w_{2})(w_{3}+w_{4)}\\
&\qquad\qquad +(v_{1}^{2}+v_{2}^{2})(w_{3}+w_{4})(w_{5}+w_{6)}]\leq 0,
\end{split}
\end{equation}
and the second line 
\begin{equation}
\begin{split}
&\frac{4v_{2}^{2}(v_{1}^{2}+v_{2}^{2})}{D^{4}}[(v_{3}^{2}+v_{4}^{2})(w_{1}+w_{2})(w_{5}+w_{6)}+(v_{5}^{2}+v_{6}^{2})(w_{1}+w_{2})(w_{3}+w_{4)}\\
&\qquad\qquad +(v_{1}^{2}+v_{2}^{2})(w_{3}+w_{4})(w_{5}+w_{6)}]\leq 0,
\end{split}
\end{equation}
and, in fact, we have only one condition

\begin{equation}\label{B12}
\begin{split}
(v_{1}^{2}+v_{2}^{2})(w_{3}+w_{4})(w_{5}+w_{6})+(v_{3}^{2}+v_{4}^{2})(w_{1}+w_{2})(w_{5}+w_{6})+(v_{5}^{2}+v_{6}^{2})(w_{1}+w_{2})(w_{3}+w_{4})\leqslant 0,
\end{split}
\end{equation}
or, more explicitly, 
\begin{equation}\label{C4}
\begin{split}
[\alpha^{2}_{1}k^{2}_{1}&y^{2k_{1}-2}_{1}+\alpha^{2}_{2}k^{2}_{2}y^{2k_{2}-2}][\alpha_{3}k_{3}(2k_{3}-1)y^{2k_{3}-2}_{3}+\alpha_{4}k_{4}(2k_{4}-1)y^{2k_{4}-2}_{4}][\alpha_{5}k_{5}(2k_{5}-1)y^{2k_{5}-2}_{5}\\
&+\alpha_{6}k_{6}(2k_{6}-1)y^{2k_{6}-2}_{6}]+[\alpha^{2}_{3}k^{2}_{3}y^{2k_{3}-2}_{3}+\alpha^{2}_{4}k^{2}_{4}y^{2k_{4}-2}_{4}]
[\alpha_{1}k_{1}(2k_{1}-1)y^{2k_{1}-2}_{1}\\
&+\alpha_{2}k_{2}(2k_{2}-1)y^{2k_{2}-2}_{2}][\alpha_{5}k_{5}(2k_{5}-1)y^{2k_{5}-2}_{5}
+\alpha_{6}k_{6}(2k_{6}-1)y^{2k_{6}-2}_{6}]+[\alpha^{2}_{5}k^{2}_{5}y^{2k_{5}-2}_{5}\\
&+\alpha^{2}_{6}k^{2}_{6}y^{2k_{6}-2}][\alpha_{1}k_{1}(2k_{1}-1)y^{2k_{1}-2}_{1}+\alpha_{2}k_{2}(2k_{4}-1)y^{2k_{4}-2}_{4}][\alpha_{3}k_{3}(2k_{3}-1)y^{2k_{3}-2}_{3}\\
&+\alpha_{4}k_{4}(2k_{4}-1)y^{2k_{4}-2}_{4}]\leq 0.
\end{split}
\end{equation}
This inequality express \ref{B10}. But, in \ref{B10}, $B_{j,k}$ and $B^{j,k}$ are functions of $(x_{1},\;x_{2})$ defined on $D\subset \R^{2}.$ In particular, $y_{i}=y_{i}(f(x))=f_{i}(x).$\\
\underline{First case}  If $(y_{1},\;y_{2})=(f_{1}(x),\;f_{2}(x))$ and $(f_{3}(x),\;f_{4}(x))$ and $(f_{5}(x),\;f_{6}(x))$ are not identically zero on $D,$ (we remark, by holomorphy, if, for example, $f_{3}(x)$ no zero on $D,$ $f_{4}(x)$ has the same property) then 
the condition \ref{C4} implies that $w_{1}+w_{2},\;w_{3}+w_{4},\;w_{5}+w_{6}$ are not all strictly positive or strictly negative. And, therefore, $\alpha_{1},\alpha_{2},\alpha_{3},\alpha_{4},\alpha_{5},$ and $\alpha_{6}$  are not all strictly positive or strictly negative. So, the hypersurface $\rho=0$ is not compact, as the theorem of Diederich and Fornaess says. Nevertheless, the necessary condition \ref{B12} is strictly stronger than the Diederich-Fornaess condition (that is to say $\{\rho =0\}$ no compact).\\
\underline{Second case} If $(f_{1}(x),\;f_{2}(x))=0\;\;on\;\; D,$ or $(f_{3}(x),\;f_{4}(x))=0\;\;on\;\; D,$ or $(f_{5}(x),\;f_{6}(x))=0\;\;on\;\; D,$ for example $(f_{5}(x),\;f_{6}(x))=0\;\;on\;\; D,$ then $f$ is, in fact, an holomorphic function from (an open set in) $\R^{2}$ to $\R^{4}.$ We can resume the preceding calculations in this case, or, more simply, consider the quadratic forms \ref{A11} when $p^{5}_{1}=p^{6}=0.$ Then, the quadratic forms \ref{A11} are only
\begin{equation}
-B_{2,2}[(p^{3}_{1})^{2}+(p^{4}_{1})^{2}]\;\; and\;\;B^{2,2}[(p^{3}_{1})^{2}+(p^{4}_{1})^{2}]
\end{equation}
which have to be neither positive definite nor negative definite. Therefore, 
\begin{equation}
B_{2,2}=B^{2,2}=0,
\end{equation}
which gives, from \ref{A15}, 
\begin{equation}\label{A17}
\begin{cases}
&v_{1}[(v_{3}^{2}+v_{4}^{2})(w_{1}+w_{2})+(v_{1}^{2}+v_{2}^{2})(w_{3}+w_{4})]=0\\
&v_{2}[(v_{3}^{2}+v_{4}^{2})(w_{1}+w_{2})+(v_{1}^{2}+v_{2}^{2})(w_{5}+w_{6})]=0.
\end{cases}
\end{equation}
If $v_{1}=0$ or $v_{2}=0,$ for example $v_{1}=0,$ then $\alpha_{1}=0,$  (so $H=\{\rho=0\}$ is not compact) or $y_{1}=f_{1}(x)=0.$\\ 
If $y_{1}=f_{1}(x)=0,$ then $y_{2}=f_{2}(x)=C$ a constant, so $(f_{3},f_{4})$ is an holomorphic function from $D\subset \R^{2}$ to $H\cap \{y_{1}=0\}\cap \{y_{2}=C\}$ which is of dimension 1 : impossible (except if $(f_{3},f_{4})$ constant). We have, of course, similar results if $f_{3}(x)=0$ or $f_{5}(x)=0.$\\
Therefore, we suppose $v_{1}\neq 0$ and the first line in \ref{A17} becomes 

\begin{equation}
(v_{3}^{2}+v_{4}^{2})(w_{1}+w_{2})+(v_{1}^{2}+v_{2}^{2})(w_{3}+w_{4})=0,
\end{equation}

that is to say 
\begin{equation}
\begin{split}
[2\alpha_{1}k_{1}y^{2k_{1}-2}_{1}(2k_{1}-1)&+2\alpha_{2}k_{2}y^{2k_{2}-2}_{2}(2k_{2}-1)](v^{2}_{3}+v^{2}_{4})\\
&+[2\alpha_{3}k_{3}y^{2k_{3}-2}_{3}(2k_{3}-1)+2\alpha_{4}k_{4}y^{2k_{4}-2}_{4}(2k_{4}-1)](v^{2}_{1}+v^{2}_{2})=0
\end{split}
\end{equation}
whereas $v^{2}_{1}+v^{2}_{2}>0$ and $v^{2}_{3}+v^{2}_{4}>0.$ So, 
\begin{equation}
\frac{2\alpha_{1}k_{1}y^{2k_{1}-2}_{1}(2k_{1}-1)+2\alpha_{2}k_{2}y^{2k_{2}-2}_{2}(2k_{2}-1)}{v_{1}^{2}+v_{2}^{2}}=-\frac{2\alpha_{3}k_{3}y^{2k_{3}-2}_{3}(2k_{3}-1)+2\alpha_{4}k_{4}y^{2k_{4}-2}_{4}(2k_{4}-1)}{v_{3}^{2}+v_{4}^{2}}, 
\end{equation}
so, $\alpha_{1},\;\alpha_{2},\;\alpha_{3},\;\alpha_{4},$ are not all strictly positive or all strictly negative. Consequently, $H=\{\rho=0\}$ is not compact, as Diederich and Fornaess \cite{DF} say.

\subsection{The successive torsion}
Unfortunately, in general the tableau $A$ is not in involution and it happens that we have to consider prolongation tableaux ( for example in non almost complex analysis) to obtain the involution of the tableau. This is always possible by corollary $3.4$ but it can appear at each step new torsion, which is not the derivatives of the previous torsion, until the first prolongation of the tableau in involution. 
Therefore, we return just after theorem \ref{C6}, and we attempt to compute the successive torsions.
Recall that we have the following  alternative: $
\mathfrak{D_{0}}=0\;\; or\;\;\mathfrak{D_{0}}\neq 0$ (see \ref{B9}).\\
\medskip

\textit{First case} : $\mathfrak{D_{0}}=0.$\\
Then, $A^{(0)}$ is involutive, and, if there is a solution to \ref{B0}, the torsion at the order 1 vanishes and therefore  
\begin{equation}\label{A4}
\begin{cases}
c^{1}_{1,2}-\sum_{i=3}^{2n}\gamma^{1}_{i}c^{i}_{1,2}=0\\
c^{2}_{1,2}-\sum_{i=3}^{2n}\gamma^{2}_{i}c^{i}_{1,2}=0.
\end{cases}
\end{equation}
A fundamental result ensures that all the successive  torsions are zero for all the orders higher than 1 if  $A^{(0)}$ is involutive.
\medskip

\textit{Second case} : $\mathfrak{D_{0}}\neq 0.$\\
Then $A^{(0)}$ is not involutive, and if $q_{0}=Rank(\mathfrak{D_{0}},\mathfrak{D_{0}}\beta,...,\mathfrak{D_{0}}\beta^{2n-3}),$ then $A^{(q)}$ is non involutive if $q<q_{0},$ and $A^{(q)}$ is involutive if $q\geq q_{0}.$ So (see \cite{BCGGG}, p. 333), the torsion of the q-prolongation of \ref{B0} is vanishing if $q>q_{0}.$ As necessary conditions to have a solution of \ref{B0}, we first have \ref{A30},
\begin{equation}
\rho_{1}\Big[c^{1}_{1,2}-\sum_{i=3}^{2n}\gamma^{1}_{i}c^{i}_{1,2}\Big]=\rho_{2}\Big[c^{2}_{1,2}-\sum_{i=3}^{2n}\gamma^{2}_{i}c^{i}_{1,2}\Big]
\end{equation}
and, also, we write the nullity of the torsion of the q-prolongation of \ref{B0}, with $q\leq q_{0}.$ Of course, this is not easy to explicit in the general case, but we can precise this. \\
We are first looking for the torsion of the first prolongation.\\
From $p_{2}=\mathcal{A}p_{1},$ we deduce 
\begin{equation}\label{A8}
\begin{cases}
p_{1,2}&=\frac{\partial \mathcal{A}}{\partial x_{1}}p_{1}+\mathcal{A}p_{1,1}\\
p_{2,2}&=\frac{\partial \mathcal{A}}{\partial x_{2}}p_{1}+\mathcal{A}p_{1,2}\\
&=\frac{\partial \mathcal{A}}{\partial x_{2}}p_{1}+\mathcal{A}\frac{\partial \mathcal{A}}{\partial x_{1}}p_{1}+\mathcal{A}^{2}p_{1,1}.
\end{cases}
\end{equation}
Deriving, in $x_{1}$ and $x_{2},$ for the second time, the first line of \ref{B2}, that is to say, deriving \ref{A6}, we obtain
\begin{equation}
\begin{cases}
\sum_{i,j=1}^{2n}\rho_{i,j}p^{1}_{1}p^{j}_{1}+\sum_{i=1}^{2n}\rho_{i}p^{i}_{1,1}=0\\
\sum_{i,j=1}^{2n}\rho_{i,j}p^{1}_{1}p^{j}_{2}+\sum_{i=1}^{2n}\rho_{i}p^{i}_{1,2}=0\\
\sum_{i,j=1}^{2n}\rho_{i,j}p^{1}_{2}p^{j}_{2}+\sum_{i=1}^{2n}\rho_{i}p^{i}_{2,2}=0,
\end{cases}
\end{equation}
and, using \ref{A8}, expressing all the derivatives of $f$ as functions of the derivatives in $x_{1},$ we have

\begin{equation}
\begin{cases}
\sum _{i=1}^{2n}\rho_{i}p^{i}_{1,1}+\sum_{i,j=1}^{2n}\rho_{i,j}p^{i}_{1}p^{j}_{1}=0\\
\sum_{i=1}^{2n}\rho_{i}\big[\sum_{k=1}^{2n}\frac{\partial \alpha_{i,k}}{\partial x_{1}}p^{k}_{1}+\alpha_{i,k}p^{k}_{1,1}\big]+\sum_{i,j=1}^{2n}\rho_{i,j}p^{i}_{1}\alpha_{j,k}p^{k}_{1}=0\\
\sum_{i,j,k,l=1}^{2n}\rho_{i,j}\big[\alpha_{i,k}p^{k}_{1}\big]\big[\alpha_{j,l}p^{l}_{1}\big]+\sum_{i=1}^{2n}\rho_{i}\big[\sum_{j=1}^{2n}\frac{\partial \alpha_{i,j}}{\partial x_{2}}p^{i}_{1}+\sum_{j,k=1}^{2n}\alpha_{i,j}\frac{\partial \alpha_{j,k}}{\partial x_{1}}p^{k}_{1}\\
\mspace{400mu} +\sum_{j,k=1}^{2n}
\alpha_{i,j}\alpha_{j,k}p^{k}_{1,1}\big]=0.
\end{cases}
\end{equation}
Using \ref{A10}, we obtain the system 
\begin{equation}
\begin{cases}
p^{1}_{1,1}\rho_{1}+p^{2}_{1,1}\rho_{2}+p^{3}_{1,1}\rho_{3}=h_{1}(p^{3}_{1},...,p^{2n}_{1},\;p^{4}_{1,1},...,\;p^{2n}_{1,1})\\
p^{1}_{1,1}\sum_{i=1}^{2n}\rho_{i}\alpha_{i,1}+p^{2}_{1,1}\sum_{i=1}^{2n}\rho_{i}\alpha_{i,2}+p^{3}_{1,1}\sum_{i=1}^{2n}\rho_{i}\alpha_{i,3}=h_{2}(p^{3}_{1},...,p^{2n}_{1},\;p^{4}_{1,1},...,\;p^{2n}_{1,1})\\
p^{1}_{1,1}\sum_{i,j=1}^{2n}\rho_{i}\alpha_{i,j}\alpha_{j,1}+p^{2}_{1,1}\sum_{i,j=1}^{2n}\rho_{i}\alpha_{i,j}\alpha_{j,2}+p^{3}_{1,1}\sum_{i,j=1}^{2n}\rho_{i}\alpha_{i,j}\alpha_{j,3}=h_{3}(p^{3}_{1},...,p^{2n}_{1},\;p^{4}_{1,1},...,\;p^{2n}_{1,1}),
\end{cases}
\end{equation}
where $h_{1},\;h_{2},\;h_{3}$ are functions in the variables $p^{3}_{1},\;p^{4}_{1},...,\;p^{2n}_{1},\;p^{4}_{1,1},\;p^{5}_{1,1},...,\;p^{2n}_{1,1}$ which are linear in $p^{4}_{1,1},\;p^{5}_{1,1},...,\;p^{2n}_{1,1}.$
So, when the determinant of this system of three equations is non zero, we can solve and obtain 
\begin{equation}
\begin{cases}
p^{1}_{1,1}=\sum_{k=4}^{2n}\gamma^{1,1,1}_{k,1,1}p^{k}_{1,1}+\sum_{k=3}^{2n}\gamma^{1,1,1}_{k,1}p^{k}_{1}+\sum_{i,k=4}^{2n}\gamma^{1,1,1}_{i1,k1}p^{i}_{1}p^{k}_{1}\\
p^{2}_{1,1}=\sum_{k=4}^{2n}\gamma^{2,1,1}_{k,1,1}p^{k}_{1,1}+\sum_{k=3}^{2n}\gamma^{2,1,1}_{k,1}p^{k}_{1}+\sum_{i,k=4}^{2n}\gamma^{2,1,1}_{i1,k1}p^{i}_{1}p^{k}_{1}\\
p^{3}_{1,1}=\sum_{k=4}^{2n}\gamma^{3,1,1}_{k,1,1}p^{k}_{1,1}+\sum_{k=3}^{2n}\gamma^{3,1,1}_{k,1}p^{k}_{1}+\sum_{i,k=4}^{2n}\gamma^{3,1,1}_{i1,k1}p^{i}_{1}p^{k}_{1}.
\end{cases}
\end{equation}
For the initial system \ref{B4}, we had the structure equations 

\begin{equation}
\begin{cases}
\theta^{1}=df_{1}-\sum_{j=3}^{2n}\gamma^{1}_{j}p^{j}_{1}dx_{1}-\sum_{j=3}^{2n}\beta_{1,j}p^{j}_{1}dx_{2}\\
\theta^{2}=df_{2}-\sum_{j=3}^{2n}\gamma^{2}_{j}p^{j}_{1}dx_{1}-\sum_{j=3}^{2n}\beta_{2,j}p^{j}_{1}dx_{2}\\
\theta^{i}=df_{i}-p^{i}_{1}dx_{1}-\sum_{j=3}^{2n}\beta_{i,j}p^{j}_{1}dx_{2}\qquad\forall i=3,...,2n,
\end{cases}
\end{equation}
on the space of variables $M=\big\{x_{1},\;x_{2},\;f_{1},...,f_{2n},\;p^{3}_{1},...,p^{2n}_{1}\big\}.$ For the first prolongation, we now have the new structure equations 

\begin{equation}\label{A12}
\begin{cases}
\theta^{3,1}&=dp^{3}_{1}-\big[\sum_{k=4}^{2n}\gamma^{3,1,1}_{k,1,1}p^{k}_{1,1}+\sum_{k=3}^{2n}\gamma^{3,1,1}_{k,1}p^{k}_{1}+\sum_{j,k=3}^{2n}p^{j}_{1}p^{k}_{1}\gamma^{3,1,1}_{j1,k1}\big]dx_{1}\\
&\qquad -\big[\sum_{k=1}^{2n}\frac{\partial \alpha_{3,k}}{\partial x_{1}}p^{k}_{1}+\alpha_{3,k}p^{k}_{1,1}\big]dx_{2}\\
\theta^{i,1}&=dp^{i}_{1}-p^{i}_{1,1}dx_{1}-\sum_{k=1}^{2n}\big(\frac{\partial \alpha_{i,k}}{\partial x_{1}}p^{k}_{1}+\alpha_{i,k}p^{k}_{1,1}\big)dx_{2}\qquad \forall i=4,...,2n, 
\end{cases}
\end{equation}
on the space of variables $M^{1}=\big\{x_{1},\;x_{2},\;f_{1},...,f_{2n},\;p^{3}_{1},...,p^{2n}_{1},\;p^{4}_{1,1},\;p^{5}_{1,1},...,p^{2n}_{1,1}\big\}.$ In fact, in these structure equations, we should have to express $p^{1}_{1},\;p^{2}_{1},\;p^{1}_{1,1},\;p^{2}_{1,1},\;p^{3}_{1,1}$ in function of the variables in $M^{1},$ as explain before. \\
From \ref{A12}, we deduce 
\begin{equation}
\begin{cases}
d\theta^{3,1}&=\sum_{k=4}^{2n}\;A^{3,1}_{(k,1,1),1}dp^{k}_{1,1}\wedge dx_{1}+\sum_{k=4}^{2n}A^{3,1}_{(k,1,1),2}dp^{k}_{1,1}\wedge dx_{2}+\sum_{j=3}^{2n}A^{3,1}_{(j,1),1}dp^{j}_{1}\wedge dx_{1}\\
&\qquad +\sum_{j=3}^{2n}A^{3,1}_{(j,1),2}dp^{j}_{1}\wedge dx_{2}+c^{3,1}_{1,2}dx_{1}\wedge dx_{2}\\
\quad\forall i& =4,...,2n,\\
d\theta^{i,1}&=-dp^{i}_{1,1}\wedge dx_{1}+\sum_{k=4}^{2n}A^{i,1}_{(k,1,1),2}dp^{k}_{1,1}\wedge dx_{2}+\sum_{j=3}^{2n}A^{i,1}_{(j,1),2}dp^{j}_{1}\wedge dx_{2}+c^{i,1}_{1,2}dx_{1}\wedge dx_{2},
\end{cases}
\end{equation}
where we do not explicit the coefficients $A^{\bullet}_{\bullet}$ which are defined in accordance with the notations of \cite{BCGGG}, p. 130. The outstanding point is 
\begin{equation}\label{A16}
A^{i,1}_{(k,1,1),1}=-\delta^{i}_{k}
\end{equation}
so (see \cite{BCGGG}, p. 138), in the calculation of 
\begin{equation}
\begin{split}
\overline{\pi}\;:\;Span&\big(\frac{\partial }{\partial p^{3}_{1}},...,\frac{\partial }{\partial p^{2n}_{1}},\frac{\partial }{\partial p^{4}_{1,1}},...,\frac{\partial }{\partial p^{2n}_{1,1}}\big)\otimes Span(x_{1},x_{2})\\
&\qquad \longrightarrow Span\big(\frac{\partial }{\partial \theta^{1}},...,\frac{\partial }{\partial \theta^{2n}},\frac{\partial }{\partial \theta^{3,1}},...,\frac{\partial }{\partial \theta^{2n,1}}\big)\otimes Span(dx_{1}\wedge dx_{2}),
\end{split}
\end{equation}
we have, if $v=\sum_{i=1}^{2}\big[\sum_{j=3}^{2n}v^{j,1}_{i}\frac{\partial }{\partial p^{j}_{1}}+\sum_{k=4}^{2n}v^{k,1,1}_{i}\frac{\partial }{\partial p^{k}_{1,1}}\big]\otimes x_{i},$
\begin{equation}
\begin{split}
\overline{\pi}(v)&=\sum_{i=1}^{2n}\big(\sum_{j=3}^{2n}A^{i}_{(j,1),1}v^{j,1}_{2}-A^{i}_{(j,1),2}v^{j,1}_{1}\big)\frac{\partial }{\partial \theta^{i}}\\
&\qquad +\big[\sum_{j=3}^{2n}\big(A^{3,1}_{(j,1),1}v^{j,1}_{2}-A^{3,1}_{(j,1),2}v^{j,1}_{1}\big)+\sum_{k=4}^{2n}\big(A^{3,1}_{(k,1,1),1}v^{k,1,1}_{2}-A^{3,1}_{(k,1,1),2}v^{k,1,1}_{1}\big)\big]\frac{\partial }{\partial \theta^{3,1}}\\
&\qquad +\sum_{i=4}^{2n}\big[ \sum_{j=3}^{2n}\big(A^{i,1}_{(j,1),1}v^{j,1}_{2}-A^{i,1}_{(j,1),2}v^{j,1}_{1}\big)+ \big(v^{i,1,1}_{2}-\sum_{k=4}^{2n}A^{i,1}_{(k,1,1),2}v^{k,1,1}_{1}\big)\big]\frac{\partial }{\partial \theta^{i,1}}.
\end{split}
\end{equation}
The torsion vanishes if there exists $v$ such that 
\begin{equation}
\overline{\pi}(v)=c=\sum_{i=1}^{2n}c^{i}_{1,2}\frac{\partial }{\partial \theta^{i}}+\sum_{j=3}^{2n}c^{j,1}_{1,2}\frac{\partial }{\partial \theta^{j,1}},
\end{equation}
which gives

\begin{equation}
\begin{cases}
\sum_{j=3}^{2n}A^{i}_{(j,1),1}v^{j,1}_{2}-A^{i}_{(j,1),2}v^{j,1}_{1}=c^{i}_{1,2}\quad \forall i=1,...,2n\\
\sum_{j=3}^{2n}\big(A^{3,1}_{(j,1),1}v^{j,1}_{2}-A^{3,1}_{(j,1),2}v^{j,1}_{1}\big)+\sum_{k=4}^{2n}\big(A^{3,1}_{(k,1,1),1}v^{k,1,1}_{2}-A^{3,1}_{(k,1,1),2}v^{k,1,1}_{1}\big)=c^{3,1}_{1,2}\\
 \sum_{j=3}^{2n}\big(A^{i,1}_{(j,1),1}v^{j,1}_{2}-A^{i,1}_{(j,1),2}v^{j,1}_{1}\big)+v^{i,1,1}_{2}-\sum_{k=4}^{2n}A^{i,1}_{(k,1,1),2}v^{k,1,1}_{1}=c^{i,1}_{1,2}\quad \forall i=4,...,2n.
\end{cases}
\end{equation}
In this system, the $2n$ first lines are exactly \ref{A14} and give the torsion of the initial system \ref{B4}. The $2n-3$ last lines give $\forall i=4,...,2n,$
\begin{equation}
v^{i,1,1}_{2}=c^{i,1}_{1,2}+\sum_{k=4}^{2n}A^{i,1}_{(k,1,1),2}v^{k,1,1}_{1}+ \sum_{j=3}^{2n}\big(A^{i,1}_{(j,1),1}v^{j,1}_{2}-A^{i,1}_{(j,1),2}v^{j,1}_{1}\big)
\end{equation}
and, transferring to the line $2n+1,$ we obtain 
\begin{equation}
\begin{split}
\sum_{k'=4}^{2n}v^{k',1,1}_{1}&\big[\sum_{k=4}^{2n}\big(A^{3,1}_{(k,1,1),1}A^{k,1}_{(k',1,1),2}\big)-A^{3,1}_{(k',1,1),2}\big]\\
&=c^{3,1}_{1,2}-\sum_{k=4}^{2n}c^{k,1}_{1,2}A^{3,1}_{(k,1,1),1}\\
&-\sum_{j=3}^{2n}\big[v^{j,1}_{2}\big(A^{3,1}_{(j,1),1}+\sum_{k=4}^{2n}A^{3,1}_{(k,1,1),1}A^{k,1}_{(j,1),1}\big)-v^{j,1}_{1}\big(A^{3,1}_{(j,1),2}+\sum_{k=4}^{2n}A^{3,1}_{(k,1,1),1}A^{k,1}_{(j,1),2}\big)\big].
\end{split}
\end{equation}
We obtain only one condition, and, if there exist $k'=4,...,2n$ such that $\sum_{k=4}^{2n}\big(A^{3,1}_{(k,1,1),1}A^{k,1}_{(k',1,1),2}\big)-A^{3,1}_{(k',1,1),2}\neq 0,$ then the torsion may be absorbed. \\
This result is true for all the prolongations. \\
We note $[k_{1},k_{2}]$ the multi-index containing $k_{1}$ times the number $1,$ and $k_{2}$ times the number $2.$  If we have the prolongation of order $k_{0},$ of the initial system, we have the space of variables 
$$M^{k_{0}}=\{x_{1},\;x_{2},\;f_{1},...,f_{2n},\,p^{3}_{1},...,p^{2n}_{1},\;p^{4}_{1,1},...,p^{2n}_{1,1},\;p^{5}_{[3,0]},...,p^{2n}_{[3,0]},...,p^{k_{0}+3}_{[k_{0}+1],0]},...,p^{2n}_{[k_{0}+1,0]}\}$$ and all the $p^{i}_{J}$ with $J$ of length $\mid J\mid \leq k_{0}+1$ are expressed as functions of the variables $p$  in $M^{k_{0}}.$\\
Now, we are looking for the prolongation of order $k_{0}+1$ of the initial system \ref{B4}, and we carry out the same calculation as, before, for the first prolongation.\\
We have derived $\rho\circ f(x)=0\quad$ $k_{0}+1$ times and obtained 
\begin{equation}
\forall j=0,...,k_{0}+1 \qquad \sum_{i=1}^{2n}\rho_{i}p^{i}_{[j,k_{0}+1-j]}=g_{j}(p_{1},p_{1,1},...,p_{[k_{0},0]}).
\end{equation}
Deriving this one time more in $x_{1}$ and $x_{2},$ we obtain a system of $k_{0}+3$ equations
\begin{equation}
\begin{cases}
\sum_{i=1}^{2n}\rho_{i}p^{i}_{[k_{0}+2,0]}=g_{k_{0}+2}(p_{1},...,p_{[k_{0}+1,0]})\\
\sum_{i=1}^{2n}\rho_{i}p^{i}_{[j,k_{0}+2-j]}=g_{j}(p_{1},...,p_{[k_{0}+1,0]})\quad \forall j=0,...,k_{0}+1,
\end{cases}
\end{equation}
and, if $h_{j}(\widehat{p})$ refers to a function of the variables $p$ belonging to $M^{k_{0}},$ this can be written 

\begin{equation}
\begin{cases}
p^{1}_{[k_{0}+2,0]}\rho_{1}+p^{2}_{[k_{0}+2,0]}\rho_{2}+...+p^{k_{0}+3}_{[k_{0}+2,0]}\rho_{k_{0}+3}=h_{0}(\widehat{p})\\
p^{1}_{[k_{0}+2,0]}\mu^{1}_{1}+p^{2}_{[k_{0}+2,0]}\mu^{1}_{2}+...+p^{k_{0}+3}_{[k_{0}+2,0]}\mu^{1}_{k_{0}+3}=h_{1}(\widehat{p})\\
...................\\
p^{1}_{[k_{0}+2,0]}\mu^{k_{0}+2}_{1}+p^{2}_{[k_{0}+2,0]}\mu^{k_{0}+2}_{2}+...+p^{k_{0}+3}_{[k_{0}+2,0]}\mu^{k_{0}+3}_{k_{0}+3}=h_{k_{0}+2}(\widehat{p})
\end{cases}
\end{equation}
where $\mu^{i'}_{j'}$ are coefficients depending on $\rho_{i}$ and $\alpha_{i,j},$ we do not explicit. When the determinant of this system of $k_{0}+3$ equations is not zero, we can solve and express $p^{1}_{[k_{0}+2,0]},...,p^{k_{0}+3}_{[k_{0}+2,0]}$ as functions of the variables $p$ in $M^{k_{0}}.$ For the prolongation of order $k_{0}+1$ of the initial system \ref{B4}, we have a new space of variables $M^{k_{0}+1}$ which contains all the variables in $M^{k_{0}}$ plus the  variables $p^{k_{0}+4}_{[k_{0}+2,0]},\; p^{k_{0}+5}_{[k_{0}+2,0]},...,\;p^{2n}_{[k_{0}+2,0]},$ and we have new structure equations relative to $\theta^{k_{0}+2,[k_{0}+1,0]},...\theta^{2n,[k_{0}+1,0]},$  
\begin{equation}
\begin{cases}
\theta^{k_{0}+3,[k_{0}+1,0]}=dp^{k_{0}+3}_{[k_{0}+1,0]}-p^{k_{0}+3}_{[k_{0}+2,0]}dx_{1}-p^{k_{0}+3}_{[k_{0}+1,1]}dx_{2}\\
\theta^{k_{0}+j,[k_{0}+1,0]}=dp^{k_{0}+j}_{[k_{0}+1,0]}-p^{k_{0}+j}_{[k_{0}+2,0]}dx_{1}-p^{k_{0}+j}_{[k_{0}+1,1]}dx_{2}\qquad if \quad j= 4,...,2n-k_{0}.
\end{cases}
\end{equation}
In fact, $p^{k_{0}+j}_{[k_{0}+1,1]}$ has to be expressed as a function of the variables $p$ in $M^{k_{0}+1},$ and $p^{k_{0}+3}_{[k_{0}+2,0]}$ is given by the Cramer system before.  But the important point is that, if $j= 4,...,2n-k_{0},$ 

\begin{equation}
d\theta^{k_{0}+j,[k_{0}+1,0]}=-dp^{k_{0}+j}_{[k_{0}+2,0]}\wedge dx_{1}+...
\end{equation}
and, therefore, $A^{k_{0}+j,[k_{0}+1,0]}_{(k_{0}+i,[k_{0}+2,0]),1}=-\delta^{i}_{j}.$ This is the same result as \ref{A16} for the first prolongation. \\
The calculations to obtain $\overline{\pi}$ are similar to the case of the first prolongation. \\

Then, we write that the torsion vanishes if there exist $v$ such that
\begin{equation}
\overline{\pi}(v)=c=\sum_{i=1}^{2n}c^{i}_{1,2}\frac{\partial }{\partial \theta^{i}}+\sum_{j=3}^{2n}c^{j,1}_{1,2}\frac{\partial }{\partial \theta^{j,1}}+...+\sum_{k=k_{0}+3}^{2n}c^{k,[k_{0}+1,0]}_{1,2}\frac{\partial }{\partial \theta^{k,[k_{0}+1,0]}}
\end{equation}
which gives 
\begin{equation}\label{A18}
\begin{cases}
c^{k_{0}+3,[k_{0}+1,0]}_{1,2}=....\\
c^{k_{0}+j,[k_{0}+1,0]}_{1,2}=v^{k_{0}+j,[k_{0}+2,0]}_{2}-\sum_{k=k_{0}+4}^{2n}A^{k_{0}+j,[k_{0}+1,0]}_{k,[(k_{0}+2,0]),2}\;v^{k_{0}+2,0]}_{1}\quad \forall j= 4,...,2n-k_{0}.
\end{cases}
\end{equation}
From this, if $j= 4,...,2n-k_{0},$, then 

\begin{equation}
v^{k_{0}+j,[k_{0}+2,0]}_{2}=c^{k_{0}+j,[k_{0}+1,0]}_{1,2}+\sum_{k=k_{0}+4}^{2n}A^{k_{0}+j,[k_{0}+1,0]}_{k,[(k_{0}+2,0]),2}\;v^{k_{0}+2,0]}_{1}
\end{equation}
and, transferring these values to the first line of \ref{A18}, we obtain only one condition. This condition is satisfied if one of the coefficients of $v^{k,[k_{0}+2,0]}_{1}$ is different of zero, and, then, the torsion may be absorbed. This situation is exactly similar with the case of the first prolongation. \\
Remember, in all the cases, that the torsion vanishes for a prolongation of order $q$ if $q>q_{0}=Rank(\mathfrak{D},\mathfrak{D}\beta,...,\mathfrak{D}\beta^{2n-3})$ because $A^{(q_{0})}$ is then involutive. \\

\section{Sufficient conditions}
A system of PDE is said involutive if there exists an ordinary integral element (see \cite{BCGGG}, p. 107) which is equivalent ( by a non trivial result) to  the tableau $A$ is in involution and the torsion vanishes. But sometimes is easier to find directly an ordinary integral element. This notion is very important because the $Cartan-K\ddot{a}hler$ theorem, in its useful version (see \cite{BCGGG}, p. 86) says that the existence of an ordinary integral element in a point implies the existence of an integral manifold tangent, at this point, to the ordinary integral element. \\
So, a sufficient condition to have a solution of \ref{B0} is to have an ordinary integral element.\\
The exterior differential system $\mathcal{I}$ is generated as a differential ideal by the sections of  the sub-bundle $I=Span(\theta^{1},...,\theta^{2n})\subset T^{\star}M.$\\
$T^{\ast}M=Span(dx_{1},dx_{2},\theta^{1},...,\theta^{2n},dp^{3}_{1},...,dp^{2n}_{1})$ and we cosider also the dual basis of $TM$\\
$TM= Span(\frac{\partial}{\partial x_{1}},\frac{\partial}{\partial x_{2}},\frac{\partial}{\partial \theta^{1}},...,\frac{\partial}{\partial \theta^{2n}},\frac{\partial}{\partial p_{1}^{3}},...,\frac{\partial}{\partial p_{1}^{2n}}).$\\
We are looking for an almost holomorphic curve verifying \ref{B0}, so we search for a ordinary integral element $E\subset T_{m}M$ with $m\in M$ and $dim(E)=2.$ So, let $(\widetilde{e_{1}},\widetilde{e_{2}})$ be a basis of $E,$ with 
\begin{equation}
\begin{split}
\widetilde{e_{1}}=\sum_{i=1}^{2}a^{1}_{i}\frac{\partial}{\partial x_{i}}+\sum_{i'=1}^{2n}b^{1}_{i'}\frac{\partial}{\partial \theta^{i'}}+\sum_{i''=3}^{2n}c^{1}_{i''}\frac{\partial}{\partial p^{i''}_{1}}\\
\widetilde{e_{2}}=\sum_{i=1}^{2}a^{2}_{i}\frac{\partial}{\partial x_{i}}+\sum_{i'=1}^{2n}b^{2}_{i'}\frac{\partial}{\partial \theta^{i'}}+\sum_{i''=3}^{2n}c^{2}_{i''}\frac{\partial}{\partial p^{i''}_{1}}.
\end{split}
\end{equation}
\\
First $E$ has to be an integral element (see \cite{BCGGG}, p. 65), that is $\forall  \varphi \in \mathcal{I},\;\;\varphi /  E=0$ i.e. $\forall i'=2,...,2n,\quad \theta^{i'} (\widetilde{e_{1}})=\theta^{i'} (\widetilde{e_{2}})=0,$ or $b^{1}_{i'}=b^{2}_{i'}=0.$ Thus,
\begin{equation}
\begin{split}
\widetilde{e_{1}}=\sum_{i=1}^{2}a^{1}_{i}\frac{\partial}{\partial x_{i}}+\sum_{i''=3}^{2n}c^{1}_{i''}\frac{\partial}{\partial p^{i''}_{1}}\\
\widetilde{e_{2}}=\sum_{i=1}^{2}a^{2}_{i}\frac{\partial}{\partial x_{i}}+\sum_{i''=3}^{2n}c^{2}_{i''}\frac{\partial}{\partial p^{i''}_{1}}.
\end{split}
\end{equation}
But, $E$ has to be an ordinary integral element (see \cite{BCGGG}, p. 73). $m$ is an ordinary zero of $\mathcal{I}\cap \Omega^{0}(M)=\emptyset$, and we need an integral flag ${0}\subset E_{1}\subset E\subset T_{m}M$ with $E_{1}$ of dimension 1 and $K\ddot{a}hler\;\; regular$ (see \cite{BCGGG} p. 68). \\
Let $\alpha \widetilde{e_{1}}+\beta \widetilde{e_{2}}=\sum_{i=1}^{2}(\alpha a^{1}_{i}+\beta a^{2}_{i})\frac{\partial}{\partial x_{i}}+\sum_{i''=3}^{2n}(\alpha c^{1}_{i''}+\beta c^{2}_{i''})\frac{\partial}{\partial p^{i''}_{1}}$ a basis of $E_{1}.$ We are looking for a 1-form $\Omega_{1}=\gamma_{1}dx_{1}+\gamma_{2}dx_{2}+\gamma_{3}\theta^{2}+...+\gamma_{2n+1}\theta^{2n}+\gamma_{2n+2}dp^{3}_{1}+...+\gamma_{4n-1}dp^{2n}_{1}$ such that $\Omega_{1}/E_{1}\neq 0$ and $E_{1}$ is an ordinary zero of $\mathcal{F}_{\Omega_{1}}(\mathcal{I})$ (see \cite{BCGGG} p. 64). \\
First, $\Omega_{1}/E_{1}=\sum_{i=1}^{2}\gamma_{i}(\alpha a_{i}^{1}+\beta a_{i}^{2})+\sum_{i''=3}^{2n}\gamma_{2n-1+i''}(\alpha c^{1}_{i''}+\beta c^{2}_{i''})\neq 0$ gives
\begin{equation}
\alpha \Big[\sum_{i=1}^{2}\gamma_{i}a^{1}_{i}+\sum_{i''=3}^{2n}\gamma_{2n-1+i''}c^{1}_{i''}\Big]+\beta \Big[\sum_{i=1}^{2}\gamma_{i}a^{2}_{i}+\sum_{i''=3}^{2n}\gamma_{2n-1+i''}c^{2}_{i''}\Big]\neq 0.
\end{equation}
We now have to define $\mathcal{F}_{\Omega_{1}}(\mathcal{I}).$ \\
If $D=Span\Big[a_{1}\frac{\partial}{\partial x_{1}}+a_{2}\frac{\partial}{\partial x^{2}}+a_{3}\frac{\partial}{\partial \theta^{2}}+...+a_{2n+1}\frac{\partial}{\partial\theta^{2n}}+a_{2n+2}\frac{\partial}{\partial p^{3}_{1}}+...+a_{4n-1}\frac{\partial}{\partial p^{2n}_{1}}\Big]$ is a straight line in $TM,$ such that $\Omega_{1}/E_{1}\neq 0,$ we write, for $\varphi=\sum_{i'=2}^{2n}\alpha'_{i'}\theta^{i'}\in \mathcal{I}\cap \Omega^{1}(M),\quad \varphi/D=\varphi_{\Omega_{1}}(D)\Omega_{1}.$ This makes sense as explained in \cite{BCGGG} p. 68. This last equality gives $\alpha'_{2}a_{3}+...+\alpha'_{2n}a_{2n+1}=\varphi_{\Omega_{1}}(D)(\gamma_{1}a_{1}+...+\gamma_{4n-1}a_{4n-1}),$ i.e.
\begin{equation}
\varphi_{\Omega_{1}}(D)=\frac{\alpha'_{2}a_{3}+...+\alpha'_{2n}a_{2n+1}}{\gamma_{1}a_{1}+...+\gamma_{4n-1}a_{4n-1}},
\end{equation}
and $\mathcal{F}_{\Omega_{1}}(\mathcal{I})$ is the space of these functions $\varphi_{\Omega_{1}}$ with $\gamma_{1}a_{1}+...+\gamma_{4n-1}a_{4n-1}\neq 0.$  Evidently, $\varphi_{\Omega_{1}}(E_{1})=0,$ but, is $E_{1}$ an ordinary zero for $\mathcal{F}_{\Omega_{1}}(\mathcal{I})\;?$\\
Let $D'$ a straight line neighbouring $E_{1}=Span\Big(\sum_{i=1}^{2}(\alpha a^{1}_{i}+\beta a^{2}_{i})\frac{\partial}{\partial x_{i}}+\sum_{i''=3}^{2n}(\alpha c^{1}_{i''}+\beta c^{2}_{i''})\frac{\partial}{\partial p^{i''}_{1}}\Big),$
\begin{equation}
D'=Span\Big(\sum_{i=1}^{2}(\alpha a^{1}_{i}+\beta a^{2}_{i}+\varepsilon_{i})\frac{\partial}{\partial x_{i}}+\sum_{i'=1}^{2n}\varepsilon'_{i'}\frac{\partial}{\partial\theta^{i'}}+\sum_{i''=3}^{2n}(\alpha c^{1}_{i''}+\beta c^{2}_{i''}+\varepsilon''_{i''})\frac{\partial}{\partial p^{i''}_{1}}\Big)
\end{equation}
with $\varepsilon,\;\varepsilon'\;\varepsilon''$ small enough.\\
We want to obtain a finite number $q$ of functions $\varphi^{1}_{\Omega_{1}}, ...,\varphi^{q}_{\Omega_{1}}\in \mathcal{F}_{\Omega_{1}}(\mathcal{I})$ such that $d\varphi^{1}_{\Omega_{1}}, ...,d\varphi^{q}_{\Omega_{1}}$ are independent and 
\begin{equation}
\Big(\forall \varphi\in \mathcal{I},\quad \varphi(D')=0\Big) \Leftrightarrow \Big(\varphi^{1}_{\Omega_{1}}(D')= ...,=\varphi^{q}_{\Omega_{1}}(D')=0\Big) .
\end{equation}
We choose $q=2n-1,$ and, for $k=1,...,2n-1,$ we define 
\begin{equation}
\varphi^{k}=c_{k}(\gamma_{1}a_{1}+...+\gamma_{4n-1}a_{4n-1})\theta^{k+1} \quad with \quad 0\neq c_{k}\in \Omega^{0}(M).
\end{equation}
Then, 
\begin{equation}
\varphi^{1}_{\Omega_{1}}(D)=c_{1}a_{3},\quad \varphi^{2}_{\Omega_{1}}(D)=c_{2}a_{4},\quad ......\quad \varphi^{2n-1}_{\Omega_{1}}(D)=c_{2n-1}a_{2n+1},
\end{equation}
and, as said before, we want 
\begin{equation}
\begin{cases}
d\varphi^{1}_{\Omega_{1}}\wedge ...\wedge d\varphi^{2n-1}_{\Omega_{1}}\neq 0\quad at\; the\; point\: E_{1},\; and,\: consequently,\; in \; a \; a\; neighbourhood,\\
\Big(\varphi^{1}_{\Omega_{1}}(D')=\varphi^{2}_{\Omega_{1}}(D')=...=\varphi^{2n-1}_{\Omega_{1}}(D')=0\Big)\Rightarrow \Big(\forall\varphi\in \mathcal{I},\quad \varphi_{\Omega_{1}}(D')=0\Big).
\end{cases}
\end{equation}
The second condition is obviously satisfied. In fact, $\forall\varphi\in\mathcal{I},\quad \varphi=\alpha_{2}\theta^{2}+...+\alpha_{2n}\theta^{2n},$ \\
so $\varphi_{\Omega_{1}}(D')=\sum_{j=2}^{2n}\alpha_{j}\theta^{j}_{\Omega_{1}}(D')=\sum_{j=2}^{2n}\frac{\alpha_{j}}{c_{j-1}(\gamma_{1}a_{1}+...+\gamma_{4n-1}a_{4n-1})}\varphi^{j-1}_{\Omega_{1}}(D')=0.$\\
For the first condition, according to \cite{BCGGG}, p. 68, we note $G_{1}(TM,\Omega_{1})$ the open set of the straight lines $D$ in $TM$ such that $\Omega_{1}/D\neq 0.$ Then $\varphi^{k}_{\Omega_{1}}\;:\; 
G_{1}(TM,\Omega_{1})\longrightarrow \mathbb{R}$ is defined before : $\varphi^{k}_{\Omega_{1}}(D)=c_{k}a_{k+2},$ so $d\varphi^{k}_{\Omega_{1}}=c_{k}da_{k+2},$ and 
\begin{equation}
d\varphi^{1}_{\Omega_{1}}\wedge ...\wedge d\varphi^{2n-1}_{\Omega_{1}}=c_{1}...c_{2n-1}da_{3}\wedge ...\wedge da_{2n+1}\neq 0.
\end{equation}

So, $E_{1}$ is an ordinary zero of $\mathcal{F}_{\Omega_{1}}(\mathcal{I}).$ Then, $E_{1}$ is said a $K\ddot{a}hler-ordinary$ point. \\
We have also to prove that $r$ is locally constant in a neighbourhood of $E_{1}$ in $V^{0}_{1}(\mathcal{I})$ (see \cite{BCGGG} p. 67 and 68), where $V^{0}_{1}(\mathcal{I})$ denotes the subspace of $K\ddot{a}hler-ordinary$ points. \\
We have to precise. 
Let $u=\sum_{i=1}^{2}u_{i}\frac{\partial}{\partial x_{i}}+\sum_{i''=3}^{2n}u_{i''}\frac{\partial}{\partial p^{i''}_{1}}$ be a basis of $E_{1},$ and let $D'$ be a straight line neighbouring $E_{1},$ \\
$D'=Span(u^{\varepsilon})=Span\Big(\sum_{i=1}^{2}(u_{i}+\varepsilon_{i})\frac{\partial}{\partial x_{i}}+\sum_{i'=1}^{2n}\varepsilon'_{i'}\frac{\partial}{\partial \theta^{i'}}+\sum_{i''=3}^{2n}(u_{i''}+\varepsilon''_{i''})\frac{\partial}{\partial p^{i''}_{1}}\Big).$ \\
If we take $\varepsilon=0,$ then $D'=E_{1}.$\\
Following \cite{BCGGG}, p. 67, we define $H(D')=\{v\in T_{m}M\;:\; \varphi ({u_{\varepsilon}},v)=0\;\;\forall \varphi\in \mathcal{I}\cap \Omega^{2}(M)\}$ the polar space of $D',$ and $r(D')=dim(H(D'))-2.$\\
We take 
\begin{equation}
\varphi=\omega_{1}\wedge \theta^{1}+...+\omega_{2n}\wedge \theta^{2n}+\alpha_{1}d\theta^{1}+...+\alpha_{2n}d\theta^{2n}\in \mathcal{I}\cap \Omega^{2}(M),
\end{equation} 
with $\omega_{i}=d^{i}_{1}dx_{1}+...+d^{i}_{4n-1}dp^{2n}_{1} \in \Omega^{1}(M)$ and $\alpha_{i}\in \Omega^{0}(M),$ and 
\begin{equation}
v=v_{1}\frac{\partial}{\partial x_{1}}+v_{2}\frac{\partial}{\partial x_{2}}+v_{3}\frac{\partial}{\partial \theta_{1}}+...+v_{2n+2}\frac{\partial}{\partial \theta_{2n}}+v_{2n+3}\frac{\partial}{\partial p_{1}^{3}}+...+v_{4n}\frac{\partial}{\partial p_{1}^{2n}}.
\end{equation}
For simplifying notations, the basis $(dx_{1},dx_{2},\theta^{1},...,\theta^{2n},dp^{3}_{1},...,dp^{2n}_{1})$ we have seen before, is now noted $(\mu_{1},\mu_{2},\mu_{3},...,\mu_{2n+2},\mu_{2n+3},...,\mu_{4n}).$ So, we have 
\begin{equation}
\varphi =\sum_{k=1}^{2n}\sum_{i=1}^{4n-1}d^{k}_{i}\mu_{i}\wedge\theta^{k}+\sum_{k=1}^{2n}\alpha_{k}d\theta^{k},
\end{equation}
and, therefore, 
\begin{equation}
\varphi(u_{\varepsilon},v)=\sum_{k=1}^{2n}\sum_{i=1}^{4n-1}d^{k}_{i}[\mu_{i}(u_{\varepsilon})\theta^{k}(v)-\mu_{i}(v)\theta^{k}(u_{\varepsilon})]+\sum_{k=1}^{2n}\alpha_{k}d\theta^{k}(u_{\varepsilon},v).
\end{equation}
But, this is $0$ for all $\varphi\in \mathcal{I},$ i.e. $\forall d^{k}_{i}$ and $\forall \alpha_{k}.$  So $\forall k=1,...,2n,\quad \forall i=1,...,4n-1,$ we have
\begin{equation}
\begin{cases}
\mu_{i}(u_{\varepsilon})\theta^{k}(v)-\mu_{i}(v)\theta^{k}(u_{\varepsilon})=0\\
d\theta^{k}(u\varepsilon,v)=0.
\end{cases}
\end{equation}
or, more explicitly, $\forall k=1,...,2n,$ 
\begin{equation}\label{B16}
\begin{cases}
(\alpha a^{1}_{i}+\beta a^{2}_{i}+\varepsilon_{i})v_{k+1}-v_{i}\varepsilon '_{k}=0 \quad \forall i=1,2\\
\varepsilon_{i'}v_{k+1}-v_{i'+1}\varepsilon'_{k}=0 \quad \forall i'=1,...,2n\\
(\alpha c^{1}_{i''}+\beta c^{2}_{i''}+\varepsilon''_{i''})v_{k+1}-v_{2n-1+i''}\varepsilon'_{k}=0\quad \forall i''=3,...,2n\\
d\theta^{2}(u_{\varepsilon},v)=0\\
d\theta^{j}(u_{\varepsilon},v)=0\quad \forall j=3,...,2n.
\end{cases}
\end{equation}
$\alpha\widetilde{e_{1}}+\beta\widetilde{e_{2}}$ is a basis of $E_{1},$ so there exists $i$ or $i''$ such that $\alpha a^{1}_{i}+\beta a^{2}_{i}\neq 0$ or $\alpha c^{1}_{i''}+\beta c^{2}_{i''}\neq 0.$ For example, without particularizing, we suppose $\alpha a^{1}_{1}+\beta a^{2}_{1}\neq 0$ (but, if $\alpha a^{1}_{2}+\beta a^{2}_{2}\neq 0\;\; or\;\;\alpha a^{1}_{3}+\beta a^{2}_{3}\neq 0\;\; or\;\;...\;\; or\;\; \alpha a^{1}_{2n}+\beta a^{2}_{2n}\neq 0\;\; or\;\; \alpha c^{1}_{1}+\beta c^{2}_{1}\neq 0\;\; or\;\; ...\;\; or\;\; \alpha c^{1}_{2n}+\beta c^{2}_{2n}\neq 0,$ then an analogous demonstration can be obtain). And we suppose $\varepsilon$ small enough to have $\alpha a^{1}_{1}+\beta a^{2}_{1}+\varepsilon_{1}\neq 0.$ Then, from the three first lines of \ref{B16}, $\forall k=1,...,2n,$ $\forall i'=1,...,2n,$ and $\forall i''=3,...,2n,$ 
\begin{equation}
v_{k+1}=\frac{v_{1}\varepsilon'_{k}}{\alpha a^{1}_{1}+\beta a^{2}_{1}+\varepsilon_{1}}=\frac{v_{2}\varepsilon'_{k}}{\alpha a^{1}_{2}+\beta a^{2}_{2}+\varepsilon_{2}}=\frac{v_{i'+1}\varepsilon'_{k}}{\varepsilon_{i'}}=\frac{v_{2n-1+i''}\varepsilon'_{k}}{\alpha c^{1}_{i''}+\beta c^{2}_{i''}+\varepsilon''_{i''}} .
\end{equation}
(In the three last terms, if the denominator is $0,$ the numerator is also $0.$)\\
 Now, we study two cases. \\
\\
\underline{First case.}\\
There exists $k_{0}\in \{1,...,2n\}$ such that $\varepsilon'_{k_{0}}\neq 0.$ Then $\forall i'=1,...,2n,$ and $\forall i''=3,...,2n,$ 
\begin{equation}
\frac{v_{1}}{\alpha a^{1}_{1}+\beta a^{2}_{1}+\varepsilon_{1}}=\frac{v_{2}}{\alpha a^{1}_{2}+\beta a^{2}_{2}+\varepsilon_{2}}=\frac{v_{i'+1}}{\varepsilon_{i'}}=\frac{v_{2n-1+i''}}{\alpha c^{1}_{i''}+\beta c^{2}_{i''}+\varepsilon''_{i''}},
\end{equation}
and therefore, 
\begin{equation}\label{B18}
\begin{cases}
v_{2}=\frac{v_{1}(\alpha a^{1}_{2}+\beta a^{2}_{2}+\varepsilon_{2})}{\alpha a^{1}_{1}+\beta a^{2}_{1}+\varepsilon_{1}}\\
v_{i'+1}=\frac{v_{1}\varepsilon'_{i'}}{\alpha a^{1}_{1}+\beta a^{2}_{1}+\varepsilon_{1}}\quad \forall i'=1,...,2n\\
v_{2n-1+i''}=\frac{v_{1}(\alpha c^{1}_{i''}+\beta c^{2}_{i''}+\varepsilon''_{i''})}{\alpha a^{1}_{1}+\beta a^{2}_{1}+\varepsilon_{1}}\quad \forall i''=3,...,2n.
\end{cases}
\end{equation}

Using \ref{B14}, we can calculate
\begin{equation}\label{S2}
\begin{split}
d\theta^{2}(u_{\varepsilon},v)&=\sum_{i=3}^{2n}\Big(\frac{\partial\gamma^{2}_{i}}{\partial x_{2}}-\frac{\partial (\alpha_{2,1}\gamma^{1}_{i}+\alpha_{2,2}\gamma^{2}_{i}+\alpha_{2,i}}{\partial x_{1}}\Big)p^{i}_{1}\Big[(\alpha a^{1}_{1}+\beta a^{2}_{1}+\varepsilon_{1})v_{2}
 -(\alpha a^{1}_{2}+\beta a^{2}_{2}+\varepsilon_{2})v_{1}\Big]\\
&\qquad\qquad -\sum_{i=3}^{2n}\gamma^{2}_{i}\Big[(\alpha c^{1}_{i}+\beta c^{2}_{i}+\varepsilon''_{i})v_{1}-(\alpha a^{1}_{1}+\beta a^{2}_{1}+\varepsilon_{1})v_{2n-1+i}\Big]\\
&\qquad\qquad -\sum_{i=3}^{2n}(\alpha_{2,1}\gamma^{1}_{i}+\alpha_{2,2}\gamma^{2}_{i}+\alpha_{2,i})\big[(\alpha c^{1}_{i}+\beta c^{2}_{i}+\varepsilon''_{i})v_{2}-(\alpha a^{1}_{2}+\beta a^{2}_{2}+\varepsilon_{2})v_{2n-1+i}\big]
\end{split}
\end{equation}
and, $\forall j=3,...,2n,$
\begin{equation}\label{S4}
\begin{split}
d\theta^{j}(u_{\varepsilon},v)&=-\Big[(\alpha c^{1}_{j}+\beta c^{2}_{j}+\varepsilon''_{j})v_{1}-(\alpha a^{1}_{1}+\beta a^{2}_{1}+\varepsilon_{1})v_{2n-1+j}\Big]\\
&\qquad\qquad -\sum_{i=3}^{2n}\frac{\partial (\alpha_{j,1}\gamma^{1}_{i}+\alpha_{j,2}\gamma^{2}_{i}+\alpha_{j,i})}{\partial x_{1}}p^{i}_{1}\big[(\alpha a^{1}_{1}+\beta a^{2}_{1}+\varepsilon_{1})v_{2}-(\alpha a^{1}_{2}+\beta a^{2}_{2}+\varepsilon_{2})v_{1}\big]\\
&\qquad\qquad -\sum_{i=3}^{2n}(\alpha_{j,1}\gamma^{1}_{i}+\alpha_{j,2}\gamma^{2}_{i}+\alpha_{j,i})\big[(\alpha c^{1}_{j}+\beta c^{2}_{j}+\varepsilon''_{j})v_{2}-(\alpha a^{1}_{2}+\beta a^{2}_{2}+\varepsilon_{2})v_{2n-1+i}\big].
\end{split}
\end{equation}
But, from \ref{B18}, we obtain 
\begin{equation}
\begin{cases}
(\alpha a^{1}_{1}+\beta a^{2}_{1}+\varepsilon_{1})v_{2}-v_{1}(\alpha a^{1}_{2}+\beta a^{2}_{2}+\varepsilon_{2})=0\\
(\alpha a^{1}_{1}+\beta a^{2}_{1}+\varepsilon_{1})v_{2n-1+i''}-v_{1}(\alpha c^{1}_{i''}+\beta c^{2}_{i''}+\varepsilon''_{i''})=0\\
(\alpha c^{1}_{i''}+\beta c^{2}_{i''}+\varepsilon''_{i''})v_{2}-v_{2n-1+i''}(\alpha a^{1}_{2}+\beta a^{2}_{2}+\varepsilon_{2})=0.
\end{cases}
\end{equation}
These quantities appear in \ref{S2} and \ref{S4} in the expressions between $[\quad ],$ so \ref{S2} and \ref{S4} are satisfied and $v=(v_{1},...,v_{4n-1})$ given by \ref{B18} are the solutions of \ref{B16}. In other words, in this case, $dim\,H(D')=1.$ \\
We recall that we want $r(D')$ locally constant, i.e. we want $dim\;H(D')$ locally constant, and therefore we want $dim\;H(D')=1$ if $\varepsilon$ is small enough. \\

\underline{Second case.}\\
Now, $\forall k=1,...,2n,\quad \varepsilon'_{k}=0.$\\
We note $\widehat{\varepsilon}=(\varepsilon_{1},\varepsilon_{2},0,...,0,\varepsilon''_{3},...,\varepsilon''_{2n}).$ Then, \ref{B16} can be written 
\begin{equation}\label{S8}
\begin{cases}
v_{3}=...=v_{2n+1}=0\\
d\theta^{2}(u_{\widehat{\varepsilon}},v)=0\\
d\theta^{j}(u_{\widehat{\varepsilon}},v)=0\quad \forall j=3,...,2n.
\end{cases}
\end{equation}
\\
and we want $dim\;H(D')=1.$ The two last equation are calculated in \ref{S2} and \ref{S4}, and, if we note 
\begin{equation}\label{S12}
\begin{cases}
X_{2}=(\alpha a^{1}_{1}+\beta a^{2}_{1}+\varepsilon_{1})v_{2}-(\alpha a^{1}_{2}+\beta a^{2}_{2}+\varepsilon_{2})v_{1}\\
X_{i}=(\alpha c^{1}_{i}+\beta c^{2}_{1}+\varepsilon''_{i})v_{1}-(\alpha a^{1}_{1}+\beta a^{2}_{1}+\varepsilon_{1})v_{2n-1+i}\quad \forall i=3,...,2n\\
X_{2n-2+i}=(\alpha c^{1}_{i}+\beta c^{2}_{1}+\varepsilon''_{i})v_{2}-(\alpha a^{1}_{2}+\beta a^{2}_{2}+\varepsilon_{2})v_{2n-1+i}\quad \forall i=3,...,2n,
\end{cases}
\end{equation}
the two last lines of \ref{S8} become
\begin{equation}\label{S10}
\begin{cases}
&\sum_{i=3}^{2n}\Big(\frac{\partial \gamma^{2}_{i}}{\partial x_{2}}-\frac{\partial (\alpha_{2,1}\gamma^{1}_{i}+\alpha_{2,2}\gamma^{2}_{i}+\alpha_{2,i})}{\partial x_{1}}\Big)p^{i}_{1}X_{2}-\sum_{i=3}^{2n}\gamma^{1}_{i}X_{i} \\
&\qquad -\sum_{i=3}^{2n}(\alpha_{2,1}\gamma^{1}_{i}+\alpha_{2,2}\gamma^{2}_{i}+\alpha_{2,i})X_{2n-2+i}=0\\
&\; \forall j=3,...,2n,\\
& \sum_{i=3}^{2n}\frac{\partial (\alpha_{j,1}\gamma^{1}_{i}+\alpha_{j,2}\gamma^{2}_{i}+\alpha_{j,i})}{\partial x_{1}}p^{i}_{1}X_{2}+X_{j} +\sum_{i=3}^{2n}(\alpha_{j,1}\gamma^{1}_{i}+\alpha_{j,2}\gamma^{2}_{i}+\alpha_{j,i})X_{2n-2+i}=0.
\end{cases}
\end{equation}
But, the variables $X_{2},...,X_{4n-2}$ are not independent. In fact, $\forall i=3,...,2n,$ 
\begin{equation}\label{S5}
(\alpha c^{1}_{i}+\beta c^{2}_{i}+\varepsilon''_{i})X_{2}+(\alpha a^{1}_{2}+\beta a^{2}_{2}+\varepsilon_{2})X_{i}-(\alpha a^{1}_{1}+\beta a^{2}_{1}+\varepsilon_{1})X_{2n-2+i}=0.
\end{equation}
\\
We now define linear applications \\
$$f_{\varepsilon}=f\; :\; \mathbb{R}^{2n} \longrightarrow \mathbb{R}^{4n-3}$$ defined by 
$f(v_{1},v_{2},v_{2n+2},...,v_{4n-1})=(X_{2},X_{3},...,X_{2n},X_{2n+1},...,X_{4n-2})$ 
where the $X_{i}$ are defined by \ref{S12}, and 
$$g_{\varepsilon}=g\; :\; \mathbb{R}^{4n-3} \longrightarrow \mathbb{R}^{2n-1}$$ derfined by $g(X_{2},X_{3},...,X_{4n-2})=(Y_{1},Y_{2},...,Y_{2n-1})$ where $Y_{1},Y_{2},...,Y_{2n-1}$ are successively defined by the first members of the equations \ref{S10}. Then, \ref{S8} can be written 
$$\forall v=(v_{1},v_{2},v_{2n+2},...,v_{4n-1})\in \mathbb{R}^{4n},\quad     g\circ f(v)=0,$$
and we want $dim\;Ker(g\circ f)=1.$\\
First, we remark that $X_{2},...,X_{4n-2}$ are linked by \ref{S5}, so $dim\;Im(f)\leqslant 2n-1.$ \\
But, it is easy to extract, from the matrix of $f,$ a $(2n-1,2n-1)$ determinant equal to $A_{1}^{2n-1}\neq 0.$ Therefore, $dim\;Im(f)=2n-1,$ and \ref{S5} are the equation which define $Im(f)$ as a subspace of $\mathbb{R}^{4n-3}.$ So, $dim(Ker(f))=1.$ \\
We now have two cases.\\
\\
$\bullet$ If $Im(f)\cap Ker(g)=\{0\},$ then $Ker (g\circ f)=Ker(f),$ and we have $dim\;Ker(g\circ f)=1.$\\
\\
$\bullet$ If $Im(f)\cap Ker(g)\neq \{0\}, $ then $Ker(g\circ f)$ is strictly greater than $Ker(f),$ and $dim\;Ker(g\circ f)>1.$ \\
\\
Finally, $E_{1}$ is $K\ddot{a}hler-regular$ if $Ker(g)\cap Im(f)=\{0\}$ when $\alpha a^{1}_{1}+\beta a^{2}_{1}\neq 0.$ \\
\\
But, as explained before, analogous demonstrations can be obtained if $\alpha a^{1}_{i}+\beta a^{2}_{i}\neq 0$ for some $i,$ or $\alpha c^{1}_{i''}+\beta c^{2}_{i''}\neq 0$ for some $i''.$ To precise this, we give, succinctly, the case $\alpha c^{1}_{3}+\beta c^{2}_{3}\neq 0.$\\
If there exists $k_{0}\in \{2,...,2n\}$ such that $\varepsilon'_{k_{0}}\neq 0,$ then, from the three first lines of \ref{B16},
\begin{equation}\label{B28}
\begin{cases}
v_{i}=\frac{v_{2n+2}(\alpha a^{1}_{i}+\beta a^{2}_{i}+\varepsilon_{i})}{\alpha c^{1}_{3}+\beta c^{2}_{3}+\varepsilon''_{3}}\quad \forall i=1,2\\
v_{i'+1}=\frac{v_{2n+2}\varepsilon'_{i'}}{\alpha c^{1}_{3}+\beta c^{2}_{3}+\varepsilon''_{3}}\quad \forall i'=2,...,2n\\
v_{2n-1+i''}=\frac{v_{2n+2}(\alpha c^{1}_{i''}+\beta c^{2}_{i''}+\varepsilon''_{i''})}{\alpha c^{1}_{3}+\beta c^{2}_{3}+\varepsilon''_{3}}\quad \forall i''=4,...,2n.
\end{cases}
\end{equation}
From this, we deduce 
\begin{equation}\label{S22}
\begin{cases}
X_{2}=(\alpha a^{1}_{1}+\beta a^{2}_{1}+\varepsilon_{1})v_{2}-(\alpha a^{1}_{2}+\beta a^{2}_{2}+\varepsilon_{2})v_{1}=0\\
X_{i}=(\alpha c^{1}_{i}+\beta c^{2}_{1}+\varepsilon''_{i})v_{1}-(\alpha a^{1}_{1}+\beta a^{2}_{1}+\varepsilon_{1})v_{2n-1+i}=0\quad \forall i=3,...,2n\\
X_{2n-2+i}=(\alpha c^{1}_{i}+\beta c^{2}_{1}+\varepsilon''_{i})v_{2}-(\alpha a^{1}_{2}+\beta a^{2}_{2}+\varepsilon_{2})v_{2n-1+i}=0\quad \forall i=3,...,2n,
\end{cases}
\end{equation}
So, \ref{B28} gives the solutions of \ref{B16}, and $dim\;H(D')=1.$\\
If, now, $\forall k=2,...,2n,\quad \varepsilon'_{k}=0,$ \ref{B16} becomes \ref{S8}. We then define $f$ and $g$ as previously and obtain the same conclusion : if $\alpha c^{1}_{3}+\beta c^{2}_{3}\neq 0,$ then $E_{1}$ is $K\ddot{a}hler-regular$ if $Ker(g)\cap Im(f)=\{0\}.$
We obtained the theorem below
\begin{theorem} $E_{1}$ is $K\ddot{a}hler-regular$ if $Ker(g)\cap Im(f)=\{0\},$ so, then, $E$ is an ordinary integral element and we have an almost holomorphic curve.
\end{theorem}
Using the Cartan-$K\ddot{a}hler$ theorem (see \cite{BCGGG} pp. 81-86), we have directly the existence of germ of a disk in the hypersurface associated to the ordinary integral element.\qed
\medskip

For ending, we can precise that $Ker(g)\cap Im(f)=\{0\}$ if and only if the system 
\begin{equation}
\begin{cases}
&\sum_{i=3}^{2n}\Big(\frac{\partial \gamma^{2}_{i}}{\partial x_{2}}-\frac{\partial (\alpha_{2,1}\gamma^{1}_{i}+\alpha_{2,2}\gamma^{2}_{i}+\alpha_{2,i})}{\partial x_{1}}\Big)p^{i}_{1}X_{2}-\sum_{i=3}^{2n}\gamma^{1}_{i}X_{i} \\
&\qquad -\sum_{i=3}^{2n}(\alpha_{2,1}\gamma^{1}_{i}+\alpha_{2,2}\gamma^{2}_{i}+\alpha_{2,i})X_{2n-2+i}=0\\
&\; \forall j=3,...,2n,\\
& \sum_{i=3}^{2n}\frac{\partial (\alpha_{j,1}\gamma^{1}_{i}+\alpha_{j,2}\gamma^{2}_{i}+\alpha_{j,i})}{\partial x_{1}}p^{i}_{1}X_{2}+X_{j} +\sum_{i=3}^{2n}(\alpha_{j,1}\gamma^{1}_{i}+\alpha_{j,2}\gamma^{2}_{i}+\alpha_{j,i})X_{2n-2+i}=0\\
&(\alpha c^{1}_{j}+\beta c^{2}_{j}+\varepsilon''_{j})X_{2}+(\alpha a^{1}_{2}+\beta a^{2}_{2}+\varepsilon_{2})X_{j}-(\alpha a^{1}_{1}+\beta a^{2}_{1}+\varepsilon_{1})X_{2n-2+j}=0
\end{cases}
\end{equation}
is a Cramer system with unknowns $X_{2},\;X_{3},...,\;X_{4n-2},$ whose the determinant is $\neq 0.$

\section{A test for finding curves in a real analytic hypersurface}

All the results in this paragraph can be certainly adapted  in a more general setting but we chose to write its in the nearly complex setting for the convenience of the reader. Furthermore, to our knowledge, the  procedure is new. 

In the section 3 we give a necessary condition to find a germ of curve in a real analytic hypersurface. Clearly this condition is far to be sufficient but we can obtain a more precise condition. With the same notations as in section 3, $f=(f_1, \cdots ,f_{2n})$ and $p^{i}_{1}=(p^{1}_{1} , \cdots ,  p^{2n}_1)$  the necessary conditions can be viewed (in a less intrinsic way) as the following : take $(p^{i}_1)$  a 1-jet satisfying the equations (with the Einstein notation) :
\begin{equation}\label{T0}
\begin{cases}
& p^{i}_2=J(p^{i}_1)\\
&\partial_i \rho(f)p^{i}_1=0  \\
&\partial _i\rho (f)J(p^{i}_1)=0,
\end{cases}
\end{equation}
roughly speaking the torsion is the obstruction to construct a 2-jet $(p^{i}_1,p^{i}_{11})$ with $(p^{i}_1)$  as above satisfying the derivatives of the previous equations  with respect to $x_1$ and $x_2$. Then $(p^{i}_{11})$ has to verify the following equations: 
\begin{equation}
\begin{cases}
&-D\rho (f)(p^{i}_{11})=D^2\rho(f)(p^{i}_{1} ,p^{i}_1) \\
&-D\rho (f)(J(p^{i}_{11}))=D\rho(DJ(p^{i}_1)(p^{i}_1))+D^2\rho(f)(p^{i}_{1} ,J(p^{i}_1))\\
&-D\rho(f)(J(p^{i}_{11}))=D\rho(DJ(p^{i}_{1})(p^{i}_1))+D^2\rho(J(p^{i}_1),p^{i}_1)\\
&D\rho(f)(p^{i}_{11})=D\rho(DJ(Jp^{i}_1)(p^{i}_1))+  D\rho (J(DJ(p^{i}_1)(p^{i}_1)))+D^2\rho(J(p^{i}_1),J(p^{i}_1)),\\
\end{cases}
\end{equation}
where $DJ(p^{i}_1)$ or $DJ(Jp^{i}_1)$ denote the matrix which the entries are the differentials of the entries of $J$ applied on the vectors $(p^{i}_1)$ or $J(p^{i}_1)$.
It is obvious that the last system has solutions if and only if 
\begin{equation}\label{T2}
D^2\rho(f)(p^{i}_{1} ,p^{i}_1) + D\rho(DJ(Jp^{i}_1)(p^{i}_1))+ D\rho (J(DJ(p^{i}_1)(p^{i}_1)))+D^2 \rho(J(p^{i}_1),J(p^{i}_1))=0.
\end{equation}
The last quantity is called the Levi form which coincides to the usual Levi form in the integrable case. Unfortunately, the nullity of the Levi form  does not give precise informations for the existence of germ of curve in $\{\rho=0\}:=H$. 
Now  the equations \ref{T0} and \ref{T2} define a real analytic set perhaps with singularities. Nevertheless, we can stratify this set by Withney' s process and the strates are  real analytic manifolds. Each strate can be defined by equations on $(f,p^{i}_1)\in T^{J}H$  which define a PDE system.  We then  derive these equations with respect to $x_1$ and $x_2$ and thus compute the torsion of the above system. More precisely, suppose that the strate is defined  by the system  of equations
\begin{equation}\label{T4}
 g_j(f,p^{i}_1)=0
 \end{equation}
 
 with $1\leq j\leq k$. To compute the torsion, we differentiate with respect to $x_1$ and $x_2$ and we want that there exists $p^{i}_{11}$ such that 

\begin{equation}\label{T3}
\begin{cases}
& \partial_{f_l} g_j(f,p^{i}_1)p^{l}_1=-\partial_{p^{l}_1}g_j(f,p^{i}_{1})P^{l}_{11}\\
& \partial_{f_l}g_j(f,p^{i}_1)(J(p^{i}_1)_l+\partial_{p^{l}_1}g_j(f,p^{i}_{1})(DJ(p^{i}_1)(p^{i}_1))_l=-\partial_{p^{l}_1}g_j(f,p^{i}_{1})J(P^{i}_{11})_l
\end{cases} 
\end{equation}

Remark that for $(f,p^{i}_1)$ fixed on the strate, the above equations are endomorphims in $p^{i}_{11}$. Therefore the previous linear system (for  $(f,p^{i}_1)$ fixed on the strate) has solutions if and only if the terms at left hand the equality are in the image of it. So these terms have to satisfy a certain number of linear equations (with $(f,p^{i}_1)$ fixed) with coefficients depending analytically of $(f,p^{i}_1)$. The linear equations depend of the rank of the system \ref{T3} and then of the nullity of some minors determinant  of the matrix defining \ref{T3} and the non nullity of others. 

\begin{rem} We can note that if the previous endomorphism is surjective for $(f,p^{i}_1)$ fixed on each strate, the PDE system defined  by \ref{T3} is free torsion on each strate of the stratification of the variety defined by the equations  \ref{T0} and \ref{T2}. In the following we will see that it is always possible to obtain a PDE system with free torsion equivalent to \ref{T0} in finite steps. 
\end{rem}

Then   we have a partition of the variety defined by \ref{T4} obtained by the nullity of some minors determinant and the non nullity of others. The elements of this partition are sub-analytic sets (even semi-analytic sets, see for example \cite{L}) and so we can stratify again each element of the partition  by smooth real analytic variety such that the closure are sub-analytic sets. On each smooth strate, we consider  the PDE system consisting of  linear equations defining the image of \ref{T3} on the strate.  Reproducing the previous arguments, we compute  the torsion of this system of PDE (on each strate of a stratification). The process will  stop in finite time when the ad-hoc endomorphism is surjective on each strate.

In fact, we have proved the following lemma:

\begin{lemma}\label{Y0}
We have a (finite) collection of PDE systems $(S_t)$ defined by $h^{t}_l(f,\partial_{1}(f_j))=0$ and $\partial _2 f=J(\partial_1 f)$ such that  $X_t:=\{h^{t}_l(f,p^{i}_1)=0\}$ is a submanifold of the manifold defined by the equations \ref{T0} which can be identified to $T^{J}H$ and the closure of each $X_t$ is a semi-analytic set.  Furthermore the systems $(S_t)$ are free torsion and equivalent to the system associated \ref{T0} on $X_t$.

\end{lemma}

\begin{rem} If we have a solution of the system associated to \ref{T0} around a point $p\in H$ then, in any neighborhood of $p$  there exists $p_0$ such that there exits $t_0$ with $S_{t_0}$ has a solution passing through $p_0$. For all $t$, the varieties $X_t$ contain $(f,0)$ for all $f$ in $H$ but in general nothing else. Moreover, suppose that for all $t$ the tableaux associated of $S_t$ are in involution and the dimension of the tableau associated to $S_t$ is locally constant around $(f,p^{i}_1)$, therefore, for all $(f,p^{i}_1)$ in $X_t$, there exists a solution of $S_t$ passing through $(f,p^{i}_1)$ and then a solution of $\ref{T0}$.

\end{rem}

If for all $t$ the tableaux associated to $S_t$ are in involution and the dimension of the tableau associated to $S_t$ are locally constant around each point of $X_t$, since $S_t$ are free torsion, the systems $S_t$ are in involution. Consider the set $M$ of points $f$ in $H$ such that there exists a non trivial disk passing through $f$. Suppose that this set is non empty and pick up a point $f_0$ in the closure of $M$. Therefore we can choose $f_n$ in $M$ satisfying: there exists a regular disk passing through $f_n$ contained in $H$ and $f_n\rightarrow f_0$. Using the last remark, we can construct $t_n$ with the following properties: $t_n\rightarrow f_0$,  there exists a regular pseudo-holomorphic disk passing through $t_n$ with derivative $v_n$ of unit norm $1$,  there exists $t_0$ such that  $(t_n, v_n)\in X_{t_0}$. Extracting a subsequence, we can assume $(t_n,v_n)$ tends to $(f_0,v_0)$. If the dimension of the tableau associated to $S_{t_0}$ is locally constant on $X_{t_0}$ and if $X_{t_0}$ is closed then there exists a holomorphic curve passing through $(f_0, v_0)$. It is not the case in general (see the next section).

Unfortunately the "tableaux" of the systems $(S_t)$ are not in involution but we have an estimate on the dimension of its.

\begin{lemma}\label{Y1}
The dimension de $A_{S_t}$ on any point of $X_t$ is less than $2n-2$.

\end{lemma}

The estimates is obvious : the dimension of the tableau of the system associated to \ref{T0} is, on any point of $T^JH$, $2n-2$ (see section 3). The structure forms of $S_t$ are the same than the structure forms of \ref{T0}, furthermore $X_t$ are submanifolds of \ref{T0} and so $J^{\perp}_{S_t}\subset J^{\perp}_{\ref{T0}}$, the two previous facts guarantee the estimate.\qed

\begin{theorem}We have a (finite) collection of PDE systems $(\lambda_t)$ defined by 
$$h^{t}_l(f,\partial_{1}(f_j),\partial^{2}_{1}(f_j),\cdots,,\partial^{k}_{1}(f_j) )=0$$ such that  
$Y_t:=\{h^{t}_l(f,\partial_{1}(f_j),\partial^{2}_{1}(f_j),\cdots,,\partial^{k}_{1}(f_j) )=0\}$ is a submanifold of the manifold defined by the derivatives of the equations \ref{T0}  and the closure of each $Y_t$ is a semi-analytic set.  Furthermore, the systems $(\lambda_t)$ are free torsion  with tableaux in involution and equivalent to the system associated \ref{T0} on $Y_t$.

\end{theorem}

Applying the lemma \ref{Y0} and omit the index t, we have a system $S_1$ with free torsion but perhaps without the tableau in involution. Thanks to the lemma \ref{Y1}, the tableau of the  $(2n-2)$ prolongation of $S_1$ are in involution but torsion can appear  by prolongation procedure (the tableau of $S_1$ was not in  involution). Using the procedure to obtain the lemma \ref{Y0}, we get a  system $S_2:=(S_1^{2n-2}, T_1)$ with free torsion but perhaps with the tableaux not in involution. The lemma \ref{Y1} guarantees that the dimension of $S_2$ is less than $2n-2$, so by induction we have a family of system $S_k$ such that the following estimates is satisfied for all $k$:
$$dim A_{S_k^{2n-2}}\leq dim A_{S_k}\leq dim A_{S_{k-1}^{2n-2}}\leq dim A _{S^{1}_{k-1}}\leq dim A_{S_{k-1}}.$$
But all these integers are less than $2n-2$ and then, there exists $k_0$ such that $dim A_{S^{1}_k}=dim A_{S_k}$.\qed

\subsection{Some examples}

In the following, we will show how concretely works the previously stratification on some examples and also the technical limits of the method.

Consider the hypersurface model $H=\{\rho(z)=Re(z_3)+\vert f(z_1,z_2)\vert^2=0\}$ in $\R^{6}$, equipped with the standard complex structure defined by $i$, and  $f$ the holomorphic function defined by $f(z_1,z_2)=z_{1}^2-z_{2}^3$. We introduce the standard complexification and we denote $w_l:=(p^{j}_1+iJ(p^{j}_1))_l$ the $l$-th component of the vector $(p^{j}_1+iJ(p^{j}_1))$ and $w^{(j)}_l$ is formally the $l$-th component of  the derivatives of $w^{(j-1)}$ with respect to $\frac{\partial}{ \partial t}:=\frac{1}{ 2}(\frac{\partial}{\partial x} -i\frac{\partial}{\partial y})$. Then The previous construction gives at the first step, the following  real analytic subset of $T^{\C}H$:

\begin{equation}\label{T5}
\begin{cases}
& \rho(z)=0\\
& w_3=0\\
&\sum_i\partial_if(z_1,z_2)w_i=0.
\end{cases} 
\end{equation}
There are two strates in the Withney's stratification: 
\begin{equation}\label{T6}
\begin{cases}
& \rho(z)=0\\
& w_3=0\\
&\sum_i\partial_if(z_1,z_2)w_i=0\\
&\vert\partial_{z,w}\big(\sum_i\partial_if(z_1,z_2)w_i\big)\vert^2 > 0
\end{cases} 
\end{equation}
and 
\begin{equation}\label{T7}
\begin{cases}
& \rho(z)=0\\
& w_3=0\\
&\sum_i\partial_if(z_1,z_2)w_i=0\\
&\partial_{z,w}\big(\sum_i\partial_if(z_1,z_2)w_i\big)=0.
\end{cases} 
\end{equation}
The system associated to  \ref{T7} can be written:

\begin{equation}
\begin{cases}
& Re(z_3)=0\\
& w_3=0\\
& z_1=0\\
& z_2=0\\
& w_1=0\\
& w_2=0,
\end{cases}
\end{equation}

and it is clear that the system associated to this manifold is free torsion and in involution so there exists  holomorphic curves contained in \ref{T7} passing through $(z,0)$ but there is no non-trivial holomorphic curves passing through $(Im(z_3),0,0)$ contained in \ref{T7}. The manifold \ref{T6} is more interesting. We want to calculate the torsion of this system and therefore we differentiate with respect $\frac{\partial} { \partial t}$ where $t$ is the variable in the unit disk: 

\begin{equation}\label{Z1}
\begin{cases}
& \rho(z)=0\\
& w_3=0\\
& w^{(1)}_3=0\\
& 2z_1w_1+3z_{2}^2w_2=0\\
&2w_{1}^2+2z_1w^{(1)}_1+6z_2w_2^{2}+3z_{2}^{2}w^{(1)}_2=0\\
&\vert z_1\vert^2+\vert z_2\vert^2+\vert w_1\vert^2 >0,
\end{cases}
\end{equation} 

this  system in $(w^{(1)}_1,w^{(1)}_2,w^{(1)}_3)$  has always solutions if $z^{'}:=(z_1,z_2)\not =(0,0)$ and therefore the PDF system associated to \ref{T6} is free torsion. Furthermore,  the tableau associated to the linear system in  $(w^{(1)}_1,w^{(1)}_2,w^{(1)}_3)$ (when $(z,w)$ are fixed and satisfied \ref{T6})  is in involution, $dim A=dim A^1=1$ if $z^{'}\not = 0$ and $dim A=dim A^1=2$ if not. Then for all $(z,w)$ with $z^{'}\not =0$ such that \ref{T6} is satisfied, there is an holomorphic curve passing through $(z,w)$ at the time zero using  the main theorem on linear Pfaff's system (see \cite{BCGGG}, pp. 140).\qed 

We can remark the following fact: the point $(z,w)$ with $z^{'} =0$ and $w_1 \not = 0$ lies on the manifold defined by \ref{T6} but at this point the dimension of $A$ is not locally constant so we cannot apply the previous theory and effectively, there does not exist non trivial holomorphic curve passing through it. So the hypothesis on the locally constant dimension on $A$ is necessary in general.

Nevertheless, there exists a non trivial holomorphic curve contained in $H$ at the point $(z,0)$ with $z^{'}=0$ and $z_3=ia$ with a real number $a$. But this point lies on the closure of the strate defined by \ref{T6} and we are not able to handle this case with the previous theory, which show the limits of the method. In this simple case, it is easy to override the trouble. If we have a curve passing through this point,  it stays on $z_3$ equal  constant (thanks to the equations \ref{Z1}). Therefore the curve stays in the set where  $\vert z_{1}^2-z_{2}^3\vert$ is constant, equal to $Re( z_3)=0$. Now, it is obvious  that the curve  is $t\rightarrow (t^3,t^2,ia)$.


\bigskip

Consider the hypersurface $H=\{\rho(z)=2Re(z_3)+\vert z_1\vert ^2-\vert z_2\vert ^2=0\}$ in $\R^{6}$. In this case, it is obvious that the first step of the previous construction gives :
\begin{equation}
\begin{cases}
& \rho(z)=0\\
& w_3+w_1\bar{z_1}-w_2\bar{z_2}=0\\
& \vert w_1\vert ^2-\vert w_2\vert ^2=0;
\end{cases}
\end{equation} 
we have two strates in the Withney's stratification:
\begin{equation}\label{T8}
\begin{cases}
& \rho(z)=0\\
& w_3+w_1\bar{z_1}-w_2\bar{z_2}=0\\
& \vert w_1\vert ^2-\vert w_2\vert ^2=0\\
& \vert w_1\vert^2+\vert w_2\vert^2 >0
\end{cases}
\end{equation} 

and $w_1=w_2=0$. Obviously, the second strate does not contain any non trivial holomorphic curves; the first case is more interesting, the first prolongation of the system is done by:

\begin{equation}\label{T9}
\begin{cases}
& \rho(z)=0\\
& w_3+w_1\bar{z_1}-w_2\bar{z_2}=0\\
& \vert w_1\vert ^2-\vert w_2\vert ^2=0\\
 & w^{(1)}_3+w^{(1)}_1\bar{z_1}-w^{(1)}_2\bar{z_2}=0\\
& w^{(1)}_1\bar{w_1}-w^{(1)}_{2}\bar{w_2}=0\\
&\vert w^{(1)}_1\vert ^2-\vert w^{(1)}_2\vert ^2=0\\
&\vert w_1\vert^2+\vert w_2\vert^2 >0.
\end{cases}
\end{equation}
The last system has always solutions in $(w^{1}_1,w^{1}_{2},w^{1}_3)$ then \ref{T8} is free torsion; but, since the dimension of the tableau associated to \ref{T8} is 3 and the dimension of the tableau associated to \ref{T9} is 2 (the equation $\vert w^{(1)}_1\vert ^2-\vert w^{(1)}_2\vert ^2=0$ is contained in the equations  $\vert w_1\vert ^2-\vert w_2\vert ^2=0$  and $ w^{(1)}_1\bar{w_1}-w^{(1)}_{2}\bar{w_2}=0$) the tableau associated to \ref{T8} is not in involution. Consider the second prolongation of the system:
 
\begin{equation}\label{T10}
\begin{cases}
& \rho(z)=0\\
& w_3+w_1\bar{z_1}-w_2\bar{z_2}=0\\
& \vert w_1\vert ^2-\vert w_2\vert ^2=0\\
& w^{(1)}_1\bar{w_1}-w^{(1)}_{2}\bar{w_2}=0\\
&\vert w^{(1)}_1\vert ^2-\vert w^{(1)}_2\vert ^2=0\\
& w^{(2)}_1\bar{w_1}-w^{(2)}_{2}\bar{w_2}=0\\
& w^{(2)}_1\bar{w^{(1)}_1}-w^{(2)}_{2}\bar{w^{(1)}_2}=0\\
&\vert w^{(2)}_1\vert ^2-\vert w^{(2)}_2\vert ^2=0\\
&\vert w_1\vert^2+\vert w_2\vert^2 >0,
\end{cases}
\end{equation}
since the equation $w^{(2)}_1\bar{w}^{(1)}_1-w^{(2)}_{2}\bar{w}^{(1)}_2=0$ is contained in the three equations $w^{(1)}_1\bar{w_1}-w^{(1)}_{2}\bar{w_2}=0$,  $ w^{(2)}_1\bar{w_1}-w^{(2)}_{2}\bar{w_2}=0$, $\vert w^{(1)}_1\vert ^2-\vert w^{(1)}_2\vert ^2=0$, and the equation $\vert w^{(2)}_1\vert ^2-\vert w^{(2)}_2\vert ^2=0$ is linked with  the equations  $\vert w_1\vert ^2-\vert w_2\vert ^2=0$ and $w^{(2)}_1\bar{w_1}-w^{(2)}_{2}\bar{w_2}=0$, the  system defined by \ref{T9} is free torsion. Furthermore, for the same reasons, the dimension of the tableau associated to \ref{T10} is 2 and so the  system associated to \ref{T9} is involutive and we have holomorphic disks passing through   all $(z,w)$ such that \ref{T8} is satisfied thanks to the main theorem on linear Pfaff's system (see \cite{BCGGG}, pp. 140).\qed


\begin{thebibliography}{<00>}




\bibitem[BCGGG]{BCGGG}
R.L. Bryant, S.S. Chern, R.B. Gardner, H.L. Goldschmidt, P.A. Griffiths, 
\newblock Exterior Differential Systems, 
\newblock Spinger Verlag (1991).

\bibitem[DA1]{DA1}
D'Angelo, John P,
\newblock Real hypersurfaces, orders of contact, and applications, 
\newblock Ann. of Math. (2) 115 (1982), no. 3, 615–637.

\bibitem[DA2]{DA2}
 D'Angelo, John P,
 \newblock Several complex variables and the geometry of real hypersurfaces, 
 \newblock  Studies in Advanced Mathematics. CRC Press, Boca Raton, FL, 1993. xiv+272 pp.
 
 \bibitem[DF]{DF}
 Diederich, Klas; Fornaess, John E, 
 \newblock Pseudoconvex domains with real-analytic boundary, 
 \newblock Ann. of Math. (2) 107 (1978), no. 2, 371–384.
 

\bibitem[L]{L}
Łojasiewicz, Stanislas,
\newblock Sur la géométrie semi- et sous- analytique,
\newblock Annales de l’institut Fourier, tome 43, no 5 (1993), p. 1575-1595.







\end{thebibliography}
\end{document}